\documentclass[11pt,twoside,a4paper]{article}
\usepackage[T1]{fontenc}

\usepackage[australian]{babel}

\usepackage[driver=pdftex,margin=3cm,heightrounded=true,centering]{geometry}

\usepackage{hyperref}

\usepackage{amssymb,amsmath,amsfonts}
\usepackage{amsthm}
\usepackage[capitalise]{cleveref}

\allowdisplaybreaks
\numberwithin{equation}{section}

\usepackage[dvipsnames]{xcolor}
\definecolor{MyBlue}{cmyk}{1,0.13,0,0.63}
\definecolor{MyGreen}{cmyk}{0.91,0,0.88,0.52}
\definecolor{MyRed}{rgb}{.6,0,0}
\newcommand{\mylinkcolor}{MyBlue}
\newcommand{\mycitecolor}{MyGreen}
\newcommand{\myurlcolor}{MyRed}

\makeatletter
\def\@endtheorem{\endtrivlist}% NEW
\makeatother
\theoremstyle{plain}
\newtheorem{thm}{Theorem}[section]
\newtheorem{lem}[thm]{Lemma}
\newtheorem{prop}[thm]{Proposition}

\theoremstyle{definition}
\newtheorem{defn}[thm]{Definition}
\newtheorem{remark}[thm]{Remark}
\newtheorem{assumption}[thm]{Assumption}
\newtheorem{example}[thm]{Example}

\renewcommand{\eqref}[1]{\labelcref{#1}}
\crefname{thm}{Theorem}{Theorems}
\crefname{lem}{Lemma}{Lemmas}
\crefname{prop}{Proposition}{Propositions}
\crefname{coro}{Corollary}{Corollaries}
\crefname{defn}{Definition}{Definitions}
\crefname{example}{Example}{Examples}
\crefname{remark}{Remark}{Remarks}

\hypersetup{%
  bookmarksnumbered=true,bookmarksopen=false,%
  plainpages=false,% necessary to prevent duplicate page identifiers
  linktocpage=true,%
  colorlinks=true,breaklinks=true,%
  linkcolor=\mylinkcolor,citecolor=\mycitecolor,urlcolor=\myurlcolor,%
  pdfpagelayout=OneColumn,%
  pageanchor=true,%
}

\makeatletter
\def\thm@space@setup{%
  \thm@preskip=4pt plus 2pt minus 2pt
  \thm@postskip=\thm@preskip
}
\renewenvironment{proof}[1][\proofname]{\par
  \pushQED{\qed}%
  \normalfont \topsep4\p@\relax % NEW
  \trivlist
  \item[\hskip\labelsep
        \itshape
    #1\@addpunct{.}]\ignorespaces
}{%
  \popQED\endtrivlist\@endpefalse
}
\makeatother

\usepackage{enumitem}
\setlist{topsep=4pt plus 2pt minus 2pt,partopsep=0pt,itemsep=2pt plus 2pt minus 2pt,parsep=0.5\parskip}
\setenumerate[1]{label=(\arabic*)}

\usepackage{fancyhdr}
\newcommand{\MR}[1]{}

\let\OLDthebibliography\thebibliography
\renewcommand\thebibliography[1]{
  \addcontentsline{toc}{section}{\refname}
  \OLDthebibliography{CPRS06b}
  \setlength{\parskip}{0pt}
  \setlength{\itemsep}{0pt plus 0.3ex}
}

\usepackage{standalone}

\usepackage{pgfplots}
\usepgfplotslibrary{fillbetween}
\pgfplotsset{compat=1.14}
\pgfplotsset{
    width = 0.99\linewidth,
    xlabel = $x$,
    ylabel = $y$,
    xmin = -1,xmax=1,
    ymin = -7,ymax=7,
    xtick = {-1,0,1},
    axis line style={draw=none},
    tick label style={font=\tiny},
    label style={font=\footnotesize},
}

\usepackage{subcaption}

\RequirePackage{slashed}

\newcommand{\N}{\mathbb{N}}
\newcommand{\R}{\mathbb{R}}
\newcommand{\C}{\mathbb{C}}
\newcommand{\Z}{\mathbb{Z}}
\newcommand{\A}{\mathcal{A}}
\newcommand{\mH}{\mathcal{H}}
\newcommand{\D}{\mathcal{D}}

\newcommand{\E}{\mathcal{E}}

\DeclareMathOperator{\Dom}{Dom}
\DeclareMathOperator{\End}{End}
\DeclareMathOperator{\Ker}{Ker}
\DeclareMathOperator{\supp}{supp}
\DeclareMathOperator{\Lip}{Lip}
\DeclareMathOperator{\ev}{ev}
\DeclareMathOperator{\Ran}{Ran}
\DeclareMathOperator{\Hom}{Hom}
\DeclareMathOperator{\Id}{Id}

\renewcommand{\Re}{\mathop{\textnormal{Re}}}

\renewcommand{\bar}[1]{\overline{#1}}
\newcommand{\Cliff}{{\mathrm{Cl}}}

\newcommand{\dvol}{\textnormal{dvol}}
\newcommand{\K}{K}
\newcommand{\KK}{K\!K}

\newcommand{\til}[1]{\widetilde{#1}}

\newcommand{\hotimes}{\mathbin{\hat\otimes}}
\newcommand{\hot}{\hotimes}
\newcommand{\la}{\langle}
\providecommand{\ra}{\rangle}%!

\newcommand{\mO}{\mathcal{O}}

\newcommand{\bigmvert}{\,\big|\,}
\newcommand{\Bigmvert}{\,\Big|\,}

\newcommand{\bundlefont}[1]{{\mathtt{#1}}}
\newcommand{\bS}{\bundlefont{S}}
\newcommand{\bE}{\bundlefont{E}}
\newcommand{\bF}{\bundlefont{F}}

\newcommand{\mattwo}[4]{
  \left(\!\!\!\begin{array}{c@{~}c}#1&#2\\ #3&#4\\\end{array}\!\!\!\right)
}

\title{The Kasparov product on submersions of open manifolds}
\author{
Koen van den Dungen%
\footnote{Email: \texttt{kdungen@uni-bonn.de}}
\\[4mm]
{\normalsize 
Mathematisches Institut}, 
{\normalsize Universit\"at Bonn}\\
{\normalsize Endenicher Allee 60, D-53115 Bonn}
}

\date{} 

\begin{document}

\maketitle

\begin{abstract}
\noindent
We study the Kasparov product on (possibly non-compact and incomplete) Riemannian manifolds. Specifically, we show on a submersion of Riemannian manifolds that the tensor sum of a regular vertically elliptic operator on the total space and an elliptic operator on the base space represents the Kasparov product of the corresponding classes in $\KK$-theory. 
This construction works in general for symmetric operators (i.e.\ without assuming self-adjointness), and extends known results for submersions with compact fibres. 
The assumption of regularity for the vertically elliptic operator is not always satisfied, but depends on the topology and geometry of the submersion, and we give explicit examples of non-regular operators. 
We apply our main result to obtain a factorisation in unbounded $\KK$-theory of the fundamental class of a Riemannian submersion, as a Kasparov product of the shriek map of the submersion and the fundamental class of the base manifold.  

\vspace{\baselineskip}

\end{abstract}

% \tableofcontents

\section{Introduction}

Consider a symmetric elliptic first-order differential operator $\D$ on a Riemannian manifold $M$. 
It was shown by Baum-Douglas-Taylor \cite{BDT89} that we obtain the corresponding $\K$-homology class $[\D] := [F_\D] \in \KK(C_0(M),\C)$ from the \emph{bounded transform} $F_\D := \D(1+\D^*\D)^{-\frac12}$. 
We emphasise here that $\D$ is \emph{not} required to be self-adjoint; in fact, any closed extension of $\D$ gives the same $\K$-homology class. 
A typical example is the Dirac operator on a (possibly incomplete) Riemannian spin$^c$ manifold $M$, representing the fundamental class of $M$. 

An alternative approach was given by Higson \cite{Hig89pre} (see also \cite[\S10.8]{Higson-Roe00}), who constructed an operator $\til F_\D$ representing the same $\K$-homology class, i.e.\ such that $[\til F_\D] = [F_\D]$. The main underlying idea of Higson's approach is that \emph{locally} there is no difference between symmetric and self-adjoint operators (the difference can be noticed only on the boundary or `near infinity'). 
To be precise, given a symmetric first-order differential operator $\D$ on $M$ and a precompact open subset $U\subset M$, there exists a self-adjoint operator $\D'$ such that $\D'|_U = \D|_U$. For instance, we can pick a cut-off function $\phi \in C_c^\infty(M)$ such that $\phi|_U = 1$, and consider $\D' = \phi^* \D \phi$. Then $\D'$ is a symmetric first-order differential operator with compact support, and therefore self-adjoint. 

Higson's construction then works as follows. 
Let $\{U_j\}$ be a locally finite cover of open precompact subsets of $M$, equipped with a partition of unity $\{\chi_j^2\}$. 
For each $j$, consider a self-adjoint first-order differential operator $\D_j$ such that $\D_j|_{U_j} = \D|_{U_j}$, and consider their bounded transforms $F_{\D_j} = \D_j (1+\D_j^2)^{-\frac12}$. 
Using the partition of unity $\{\chi_j^2\}$, we construct
\begin{align*}
\til F_\D := \sum_j \chi_j F_{\D_j} \chi_j .
\end{align*}
The operator $\til F_\D$ is well-defined as a strongly convergent series, and represents the same $\K$-homology class as $F_\D$ (see \cref{thm:local_rep_KK}). We will refer to $\til F_\D$ as the \emph{localised representative} for the class $[\D]$. 
Higson then used the above construction to prove that the \emph{external} Kasparov product \cite{Kas80b} of two symmetric elliptic first-order differential operators $\D_1$ on $M_1$ and $\D_2$ on $M_2$ is represented by their tensor sum:
\begin{align*}
\begin{array}{ccccc}
\KK(C_0(M_1),\C) & \hot & \KK(C_0(M_2),\C) & \to & \KK(C_0(M_1\times M_2),\C) , \\ {} 
[\D_1] & \hot & [\D_2] & = & [\D_1\hot 1 + 1\hot\D_2] .
\end{array}
\end{align*}

In this article we will generalise Higson's result to \emph{internal} Kasparov products on submersions of open manifolds. 
A constructive approach to the internal Kasparov product in unbounded $\KK$-theory has been developed in \cite{Mes14,KL13,BMS16,MR16}. One of the pillars under this construction is Kucerovsky's theorem \cite{Kuc97}, which provides sufficient conditions allowing one to check whether an unbounded Kasparov module represents the Kasparov product. 

We consider a submersion $\pi\colon M\to B$ of (possibly non-compact) smooth manifolds $M$ and $B$. Let $\D_V$ be a \emph{vertically elliptic}, symmetric, first-order differential operator $\D_V$ on the total space $M$, and let $\D_B$ be an elliptic, symmetric, first-order differential operator on the base space $B$. 
Our main goal in this article is to construct the (internal) Kasparov product of $\D_V$ with $\D_B$. 
However, since $M$ and $B$ can be non-compact, the operators $\D_V$ and $\D_B$ may not be self-adjoint, and therefore the above-mentioned constructive approach to the internal Kasparov product does not apply. 
Nevertheless, we will show that the Kasparov product of $\D_V$ with $\D_B$ is represented by a tensor sum 
\begin{align}
\label{eq:tensor_sum}
\D := \D_V\hot 1 + 1\hot_\nabla\D_B .
\end{align}
The construction of such a Kasparov product has already been considered by Kaad and Van Suijlekom \cite{KS18,KS17bpre}, but only in the special case where the submersion $\pi\colon M\to B$ is \emph{proper}, which means that the fibres are all compact. 
The technical pillar underneath the results of \cite{KS17bpre} is a variant of Kucerovsky's theorem for half-closed modules, which was proven by Kaad and Van Suijlekom in \cite{KS19}. The proof of this theorem relies on the technical setting of ``modular cycles'' developed by Kaad \cite{Kaa15pre}. One significant drawback of the theorem from \cite{KS19} is that it only applies to the case where $\D_V$ is self-adjoint, and to ensure self-adjointness of $\D_V$ one needs to assume (for instance) that the fibres of the submersion are compact (as is done in \cite{KS17bpre}). 

In this article, we instead take the rather natural approach of extending Higson's construction of a localised representative to vertical operators. 
Our methods are therefore different and independent from those developed in \cite{KS19}. 
Using such a localised representative, we will prove that the above tensor sum $\D$ indeed represents the Kasparov product of $\D_V$ with $\D_B$. Our main result thus not only provides a new, independent proof of the results of \cite{KS18,KS17bpre}, but more importantly also extends these results to submersions with \emph{non-compact} fibres. 

Let us provide an outline of this article. 
In \cref{sec:submersion}, we consider a submersion $\pi\colon M\to B$ of smooth manifolds $M$ and $B$, and a \emph{vertically elliptic}, symmetric, first-order differential operator $\D_V$ on $M$. 
The operator $\D_V$ determines a family of operators $\{\D_b\}_{b\in B}$ acting on the fibres $M_b := \pi^{-1}(b)$ of the submersion, and correspondingly we can view $\D_V$ as an operator on a Hilbert $C_0(B)$-module. 
The first question which arises is whether $\D_V$ is a \emph{regular} operator (which is necessary for defining the operator $(1+\D_V^*\D_V)^{-\frac12}$ and the bounded transform $F_{\D_V}$). 
We will see that $\D_V$ is regular if and only if the restriction of $\Dom\D_V^*$ to $M_b$ yields a core for $\D_b^*$, for every $b\in B$. 
This property is not always satisfied, but depends on the topology and geometry of the submersion. 
In \cref{sec:regular_vertical} we give a few sufficient conditions which ensure that $\D_V$ is regular, and we provide several basic examples. 
Interestingly, we show (by example) that it is possible to obtain a regular operator even if the topology and/or geometry of the fibres changes drastically. For instance, the fibres can change from complete to incomplete manifolds, or from connected to disconnected manifolds. 
We also give several examples in which the operator $\D_V$ is \emph{not} regular, illustrating what could go wrong. 
As a general rule, it appears that changing the topology (of the fibres) by removing a submanifold of codimension $1$ only yields a regular operator if the metric is chosen such that the removed submanifold lies `at geometric infinity'. 
However, removing a submanifold of codimension $2$ (or higher) yields a regular operator for any metric. 
We leave it as an open problem to give a general characterisation of the regularity of the operator in terms of the properties of the submersion. 

A vertically elliptic, symmetric, first-order differential operator $\D_V$ on a submersion $\pi\colon M\to B$ represents a class in bivariant $\K$-theory (or $\KK$-theory) \cite{Kas80b}. However, since $\D_V$ is (in general) not self-adjoint, it does not yield an unbounded Kasparov module (or unbounded $\KK$-cycle, as defined in \cite{BJ83}). Instead, we obtain the more general notion of a half-closed module as defined by Hilsum \cite{Hil10}. 
Therefore, in order to deal with symmetric unbounded operators in bivariant $\K$-theory, we briefly review half-closed modules in \cref{sec:half-closed}. 
Here we will also prove that two half-closed modules which are `locally bounded' perturbations of each other represent the same class in $\KK$-theory (see \cref{prop:loc_bdd}). 

Subsequently, we show in \cref{sec:half-closed_submersion} that $\D_V$ indeed defines a half-closed module (provided that the regularity condition is satisfied), so that the bounded transform $F_V := \D_V(1+\D_V^*\D_V)^{-\frac12}$ of $\D_V$ defines a class $[\D_V] = [F_V] \in \KK(C_0(M),C_0(B))$. 
Next, in \cref{sec:local_representative} we generalise Higson's construction, and define a localised representative $\til F_V$ from the vertical differential operator $\D_V$. 
We prove that this localised representative $\til F_V$ represents the $\KK$-class $[\D_V]$. 

Our main goal in this article is to prove that the Kasparov product of $\D_V$ with an elliptic symmetric operator $\D_B$ on the base space $B$ is represented by the tensor sum  of \cref{eq:tensor_sum}. 
The proof relies on checking the connection and positivity conditions in the well-known theorem by Connes and Skandalis (see \cref{thm:Connes-Skandalis}). 
While the connection condition can be checked for symmetric operators without too much difficulty, the positivity condition is more problematic. 
This is where the construction of a localised representative $\til F_V$ comes into play. 
The technical heart of the proof, contained in \cref{sec:local_positivity}, consists of showing that the positivity condition is satisfied `locally'. 
Here, our approach using localised representatives offers two distinct advantages. 
First, the construction of the localised representative allows us to work locally with self-adjoint (rather than only symmetric) operators. 
Second, the construction allows to rescale each localised term independently (which is crucial in order to obtain a uniform constant in the local positivity condition). 
The construction of $\til F_V$ using a partition of unity then allows us to prove that the positivity condition is in fact satisfied globally. 
We thus prove in \cref{sec:Kasp_prod} our main result, stating that the operator $\D$ indeed represents the Kasparov product of $\D_V$ and $\D_B$. 

We apply our main result in \cref{sec:factorisation} to obtain a factorisation in unbounded $\KK$-theory of the fundamental class of a Riemannian submersion $\pi \colon M\to B$ of even-dimen\-sional spin$^c$ manifolds. The fundamental classes $[M]$ and $[B]$ (in the $\K$-homology groups $\KK(C_0(M),\C)$ and $\KK(C_0(B),\C)$) are represented by the Dirac operators $\D_M$ and $\D_B$, respectively. In $\KK$-theory we have the factorisation \cite{Con82,CS84}
\[
[M] = \pi! \otimes_{C_0(B)} [B] ,
\]
where $\pi! \in \KK(C_0(M),C_0(B))$ is the shriek map of the submersion, which can be represented by a \emph{vertical} Dirac operator $\D_V$. 
We show that this factorisation can be implemented in unbounded $\KK$-theory, meaning that $\D_M$ is unitarily equivalent to the tensor sum of $\D_V$ and $\D_B$, up to an explicit curvature term. 
This generalises the work of Kaad and Van Suijlekom \cite{KS17bpre}, who proved the same result under the additional asssumption that the submersion $\pi\colon M\to B$ is \emph{proper}, which means that the fibres $M_b = \pi^{-1}(b)$ are all compact. 
While some of our arguments closely follow the work of Kaad and Van Suij\-le\-kom, the main difference is our approach to `local positivity' (\cref{sec:local_positivity}), where the use of a localised representative (as described above) allows us to consider submersions with non-compact fibres as well. 

In this article, we only consider a special case of the Kasparov product of half-closed modules for commutative $C^*$-algebras. 
Nonetheless, this special case has not yet been dealt with in the available literature so far. 
Significant parts of our proofs remain valid for noncommutative $C^*$-algebras as well, and we intend to address the Kasparov product of half-closed modules over arbitrary (noncommutative) $C^*$-algebras in a forthcoming paper.

\subsection{Notation}

Let $E$ be a $\Z_2$-graded Hilbert module over a $\sigma$-unital $C^*$-algebra $B$. 
We denote the set of adjointable operators on $E$ as $\End_B(E)$, and the subset of compact endomorphisms as $\End_B^0(E)$. 
For any operator $T$ on $E$, we write $\deg T=0$ if $T$ is even, and $\deg T=1$ if $T$ is odd. 
The graded commutator $[\cdot,\cdot]_\pm$ is defined (on homogeneous operators) by $[S,T]_\pm := ST - (-1)^{\deg S\cdot\deg T} TS$. For $R,S,T\in\End_B(E)$, the graded commutator satisfies the following identities:
\begin{align*}
[S,T]_\pm &= - (-1)^{\deg S\cdot\deg T} [T,S]_\pm , & 
[RS,T]_\pm &= R [S,T]_\pm + (-1)^{\deg S\cdot\deg T} [R,T]_\pm S .
\end{align*}
The ordinary commutator is denoted $[\cdot,\cdot] := [\cdot,\cdot]_-$. 

For any $S,T\in\End_B(E)$ we will write $S \sim T$ if $S-T\in\End_B^0(E)$. Similarly, for self-adjoint $S,T$ we will write $S\gtrsim T$ if $S-T \sim P$ for some positive $P\in\End_B(E)$; in this case we will say that $S-T$ is \emph{positive modulo compact operators}. 

Given any regular operator $\D$, we define the \emph{bounded transform} $F_\D := \D (1+\D^*\D)^{-\frac12}$.

\subsection{Acknowledgements}

The work on this article was initiated during a short visit to the Radboud University Nij\-me\-gen (late November -- early December 2017), which was funded by the COST Action MP1405 QSPACE, supported by COST (European Cooperation in Science and Technology). 
The author thanks Walter van Suijlekom for his hospitality during this visit, and for interesting discussions and helpful comments. 
These results were first presented at the thematic program \emph{Bivariant K-theory in Geometry and Physics}, Vienna, November 2018, and the author thanks the organisers, as well as 
% the Erwin Schr\"odinger International Institute for Mathematics and Physics (ESI), 
the Erwin Schr\"odinger Institute (ESI), 
for their hospitality. 
Thanks also to Bram Mesland and Matthias Lesch for many interesting discussions. 
Finally, thanks to the referee for his/her comments.

\section{The KK-class on a submersion}
\label{sec:submersion}

We consider a smooth surjective map $\pi\colon M\to B$ between manifolds $M$ and $B$. These manifolds are allowed to be non-compact (but they are not allowed to have a boundary). 
The map $\pi\colon M\to B$ is called a \emph{submersion} if its differential $d\pi(x) \colon (TM)_x\to (TB)_{\pi(x)}$ is surjective for all $x\in M$. 
We will refer to $T_VM := \Ker d\pi$ as the \emph{vertical} tangent bundle of $M$. 

Consider a smooth (complex) vector bundle $\bE\to M$, and let $\E := \Gamma^\infty(M,\bE)$ be the $C^\infty(M)$-module of smooth sections. 
We can view $\E$ as a $C^\infty(M)$-$C^\infty(B)$-bimodule, where the right action of $f\in C^\infty(B)$ on $\psi\in\E$ is given by $(\psi f)(x) := \psi(x) f(\pi(x))$, for $x\in M$. 
Consider a first-order differential operator $\D \colon \E \to \E$. 
Its principal symbol $\sigma_\D \colon \Omega^1(M) \to \End_{C^\infty(M)}(\E)$ is given by $\sigma_\D(df) = [\D,f]$, for $f\in C^\infty(M)$. 

\begin{defn}
A first-order differential operator $\D$ on $\E$ is called \emph{vertical} if $\D$ is $C^\infty(B)$-linear. 
For a subset $U\subset M$, we say that $\D$ is \emph{vertically elliptic on $U$} if $\D$ is vertical and if for each $x\in U$ the symbol $\sigma_\D(\xi)(x) \colon \bE_x\to\bE_x$ is invertible for any non-zero $\xi \in (T_V^*M)_x$. 
If $\D$ is vertically elliptic on all of $M$, we simply say that $\D$ is vertically elliptic. 
\end{defn}

\subsection{Regularity of vertical operators}
\label{sec:regular_vertical}

Now consider a submersion $\pi\colon M\to B$ of smooth manifolds equipped with a smooth vertical metric $g_V$ (i.e.\ a hermitian structure on the real vector bundle $T_VM$). 
Then for each $b\in B$, the fibre $M_b := \pi^{-1}(b)$ carries a Riemannian metric $g_b = g_V|_{M_b}$ obtained by identifying $TM_b = T_VM|_{M_b}$, and in particular we have a corresponding volume form $\dvol_b$ on $M_b$. 
We emphasise that these Riemannian metrics $g_b$ are not assumed to be complete. 

Further, consider a smooth (complex) vector bundle $\bE\to M$ with a hermitian structure $\la\cdot|\cdot\ra_\bE$. 
We write $\bE_b := \bE|_{M_b}$ and $E_b := L^2(M_b,\bE_b)$, and we let $E_\bullet$ denote the corresponding bundle of Hilbert spaces over $B$. 
We obtain a $C^\infty(B)$-valued inner product $\la\cdot|\cdot\ra_\E$ on the $C^\infty(M)$-$C^\infty(B)$-bimodule $\E := \Gamma^\infty(M,\bE)$ by integrating along the fibres:
\[
\la\phi|\psi\ra_\E(b) := \int_{M_b} \la\phi(y)|\psi(y)\ra_\bE \dvol_b(y) ,
\]
for $\phi,\psi\in\E$. 
Then the completion of $\Gamma_c^\infty(M,\bE)$ with respect to $\la\cdot|\cdot\ra_\E$ yields the Hilbert $C_0(B)$-module $\Gamma_0(B,E_\bullet)$ of continuous sections of $E_\bullet$ vanishing at infinity. 

Let $\D$ be a vertically elliptic, symmetric, first-order differential operator on $\bE\to M$. 
We will view the closure of $\D$ as an operator on $\Gamma_0(B,E_\bullet)$ with the domain $\Dom\D= \bar{\Gamma_c^\infty(M,\bE)}$ (where the closure is taken with respect to the graph norm of $\D$). 

For each $b\in B$, consider the evaluation map $\ev_b \colon \Gamma_0(B,E_\bullet)\to E_b$, which maps $\Gamma_c^\infty(M,\bE)$ to $\Gamma_c^\infty(M_b,\bE_b)$. 
The vertical operator $\D$ restricts to a symmetric elliptic first-order differential operator $\D_b$ on $\Gamma_c^\infty(M_b,\bE_b)$. 
We consider the closure of $\D_b$ (i.e.\ the \emph{minimal extension} of $\D_b$) as an operator on the Hilbert space $E_b = L^2(M_b,\bE_b)$, and (with some abuse of notation) we denote the closure simply by $\D_b$ as well. 
Since $\ev_b(\Gamma_c^\infty(M,\bE)) = \Gamma_c^\infty(M_b,\bE_b)$, we see immediately that $\ev_b(\Dom\D)$ contains a core for $\D_b$. 
To ensure that $\D$ is a \emph{regular} operator, we will need to ensure that also $\ev_b(\Dom\D^*)$ contains a core for the adjoint of $\D_b$ (i.e. for the \emph{maximal} extension). 
To prove this, we first recall the following lemma. 

\begin{lem}[{\cite[Proposition 2.9]{Hil89}}]
\label{lem:fam_regular}
Let $\{T_b\}_{b\in B}$ be a family of self-adjoint operators on a continuous field of Hilbert spaces $\{\mH_b\}_{b\in B}$, and let $T$ be the operator on the corresponding Hilbert $C_0(B)$-module $E = \Gamma_0(B,\mH_\bullet)$ given by $(T\psi)(b) := T_b \psi(b)$ for all $\psi$ in the domain
\[
\Dom T := \{ \psi\in E : \psi(b)\in\Dom T_b \; \& \; T\psi\in E \} .
\]
If $T$ is densely defined, then $T$ is self-adjoint, and furthermore $T$ is regular if and only if $\ev_b(\Dom T)$ is a core for $T_b$ for every $b\in B$. 
\end{lem}

\begin{lem}
\label{lem:vert_regular}
The following statements are equivalent: 
\begin{enumerate}
\item the (closure of the) operator $\D$ is a regular operator on $\Gamma_0(B,E_\bullet)$;
\item for each $b\in B$, the subspace $\ev_b(\Dom\D^*) \subset \Dom\D_b^*$ is a core for $\D_b^*$;
\item for each $b\in B$, the subspace $\ev_b(\Dom\D^*)$ is equal to $\Dom\D_b^*$.
\end{enumerate}
\end{lem}
\begin{proof}
The equivalence of (1) and (2) follows by applying \cref{lem:fam_regular} to $T = \mattwo{0}{\D^*}{\D}{0}$, 
and (2) and (3) are equivalent because $\ev_b(\Dom\D^*)$ is always closed in the graph norm of $\D_b^*$. 
\end{proof}

We will describe a few special cases for which we can \emph{prove} that $\ev_b(\Dom\D^*)$ is a core for $\D_b^*$. First of all, we note that it is sufficient for every $\D_b$ to be essentially self-adjoint on $\Gamma_c^\infty(M_b,\bE_b)$ (i.e.\ the minimal extension of $\D_b$ is self-adjoint). 
This situation occurs for instance if the submersion $\pi\colon M\to B$ is proper, so that every fibre $M_b$ is compact (this is the setting studied in \cite{KS17bpre}). More generally, the self-adjointness of $\D_b$ also follows if $\D_b$ has bounded propagation speed and $M_b$ is complete. 

\begin{lem}
\label{lem:complete_regular}
If for each $b\in B$, the metric on $M_b$ is complete and $\D_b$ has bounded propagation speed, then $\D$ is regular and self-adjoint. 
\end{lem}
\begin{proof}
For each $b\in B$, the vertical operator $\D$ on $\Gamma_c^\infty(M,\bE)$ restricts to a symmetric first-order differential operator $\D_b$ on $\Gamma_c^\infty(M_b,\bE_b)$. By \cite[Proposition 10.2.11]{Higson-Roe00}, the assumptions imply that each $\D_b$ is essentially self-adjoint, and therefore $\D$ is self-adjoint. 
Moreover, we have the inclusions $\Gamma_c^\infty(M_b,\bE_b) \subset \ev_b(\Dom\D) \subset \Dom\D_b$, so in particular $\ev_b(\Dom\D^*) = \ev_b(\Dom\D)$ contains a core for $\D_b^* = \D_b$. By \cref{lem:vert_regular}, $\D$ is regular. 
\end{proof}

The assumption of fibrewise completeness of the metric is particularly suited to, for instance, Dirac-type operators (which always have bounded propagation speed). For arbitrary operators without bounded propagation speed, we may consider the following statement. 

\begin{lem}[{cf.\ \cite[2.28]{Ebe16pre}}]
\label{lem:coercive_regular}
Suppose there exists a function $f\in C^\infty(M)$ with the following properties:
\begin{itemize}
\item the map $(\pi,f)\colon M\to B\times\R$ is proper;
\item for each $b\in B$, the commutator $[\D_b,f|_{M_b}]$ is bounded on $L^2(M_b,\bE_b)$.
\end{itemize}
Then $\D$ is regular and self-adjoint. 
\end{lem}
\begin{proof}
Since $(\pi,f)\colon M\to B\times\R$ is proper, it follows that the restriction $f|_{M_b} \colon M_b\to\R$ is also proper. 
By \cite[Proposition 10.2.10]{Higson-Roe00}, the symmetric first-order differential operator $\D_b$ on $\Gamma_c^\infty(M_b,\bE_b)$ is in fact essentially self-adjoint. As in the proof of \cref{lem:complete_regular}, it follows that $\D$ is regular and self-adjoint. 
\end{proof}

Finally, we consider a case in which $\D$ is regular but in general \emph{not} self-adjoint, using a (rather strong) assumption of local triviality of $M$, $\bE$, and $\D$ over $B$. 

\begin{lem}
\label{lem:loc_triv_regular}
Suppose that the submersion $\pi\colon M\to B$, the bundle $\bE\to M$ and the operator $\D$ are \emph{locally trivial} in the following sense: for each point $b\in B$ there exist an open neighbourhood $U\subset B$ of $b$, a Riemannian manifold $N$, a hermitian vector bundle $\bF\to N$, an elliptic symmetric first-order differential operator $\D_\bF$ on $\bF\to N$, an isometry $\phi\colon\pi^{-1}(U)\to U\times N$, and a bundle isomorphism $\Phi\colon\bE|_{\pi^{-1}(U)}\to\til\bF$ covering $\phi$ (where $\til\bF$ is the pullback of $\bF\to N$ to $U\times N$) such that $\Phi \D \Phi^{-1} = \til\D_\bF$ on $\Gamma_c^\infty(U\times N,\til\bF)$ (where $\til\D_\bF$ is the pullback of $\D_\bF$). 
Then $\D$ is regular. 
\end{lem}
\begin{proof}
Consider a point $b\in B$ and an element $\eta\in\Dom\D_b^*$. By assumption, there exists a local trivialisation over a neighbourhood $U$ of $b$ (as described above). 
Pick a function $\chi\in C_c^\infty(U)$ such that $\chi(b)=1$. 
Let $\til\eta$ be the section of $\til\bF\to U\times N$ obtained as the pullback of $\Phi(\eta) \in L^2(N,\bF)$. 
Then for $x\in U$, $\psi(x) := \chi(x) \Phi^{-1}(\til\eta)$ defines an element $\psi\in\Dom\D^*$ such that $\psi(b) = \eta$. 
Hence we have shown that $\ev_b(\Dom\D^*) = \Dom\D_b^*$, and it follows from \cref{lem:vert_regular} that $\D$ is regular. 
\end{proof}

Although the above three lemmas give sufficient conditions for the regularity of $\D$, these conditions are certainly not necessary. 
In the following, we will consider a simple setup on which we can discuss several examples of both regular and non-regular operators. 

\subsubsection{Examples of (non-)regular operators}
\label{sec:regular_examples}

Let $M$ be an open subset of $(-1,1)\times\R$. We have a natural map $\pi\colon M\to (-1,1)$ given by $\pi(x,y) := x$. 
We equip $M$ with a vertical metric of the form
\[
g_V(x,y) = h(x,y) dy^2 ,
\]
where $h$ is a smooth, strictly positive function on $M$. We assume that for each $x\in(-1,1)$, $\pi^{-1}(x)$ is not empty, so that $\pi$ is a submersion. We note that each fibre $M_x = \pi^{-1}(x)$ is equipped with the Riemannian measure $\sqrt{h(x,y)} dy$. 

We consider the vertical Dirac operator $\D := -i \sqrt{h(x,y)}^{-1} \partial_y$ on $C_c^\infty(M)$. 
On each fibre $M_x$ for $x\in(-1,1)$, we obtain the operator $\D_x := -i\sqrt{h(x,y)}^{-1} \partial_y$, acting on the Hilbert space $E_x := L^2(M_x)$ with initial domain $C_c^\infty(M_x)$. The inner product on $E_x$ is given by 
\[
\la\phi|\psi\ra(x) = \int_{M_x} \bar{\phi(x,y)} \psi(x,y) \sqrt{h(x,y)} dy .
\]
Since $\D_x$ is the Dirac operator on $M_x$, we know that $\D_x$ is symmetric with respect to this inner product (which can be easily checked by a direct computation). 
We note that the domain of the adjoint is given by $\Dom\D_x^* = \{\eta\in L^2(M_x) : \sqrt{h(x,y)}^{-1} \partial_y\eta\in L^2(M_x)\}$. 
In the special case where $h=1$, we note that $\Dom\D_x^*$ equals the first Sobolev space $H^1(M_x)$, which (since $M_x$ is one-dimensional) consists of (absolutely) continuous functions. 

We view (the closure of) $\D$ as a symmetric operator on the Hilbert $C_0((-1,1))$-module $\Gamma_0((-1,1),E_\bullet)$. By \cref{lem:vert_regular}, to see if $\D$ is regular, we need to check if $\ev_x(\Dom\D^*)$ is equal to $\Dom\D_x^*$. 
In the following, we will consider several examples of $(M,g_V)$, for which we will check explicitly whether or not $\D$ is regular. 

\begin{example}[The open square]
\label{ex:open_square}
Consider the manifold $M := (-1,1)\times(-1,1)$, 
and suppose that for each $x\in(-1,1)$ there exists a constant $c_x>1$ such that $c_x^{-1} \leq h(x,y) \leq c_x$ for all $y\in(-1,1)$. 
In this case, $\D$ is regular. 
\end{example}
\begin{proof}
The assumption on $h$ ensures that $\Dom\D_x^* = H^1(-1,1)$ for each $x\in(-1,1)$. 
Consider a point $x\in(-1,1)$, and let $\eta\in\Dom\D_x^*$. 
Define $\psi(x') := \phi(x') \eta$ for all $x'\in(-1,1)$, where we have picked a function $\phi\in C_c^\infty(-1,1)$ such that $\phi(x)=1$. 
Then $\psi\in\Dom\D^*$ such that $\psi(x) = \eta$. Hence we have shown that $\ev_x(\Dom\D^*) = \Dom\D_x^*$ for any point $x\in(-1,1)$. 
\end{proof}

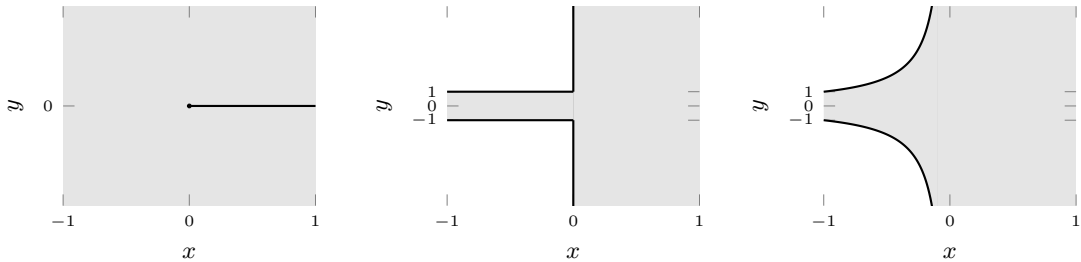
\begin{figure}[b]
\begin{subfigure}{.33\textwidth}
  \centering
  \begin{tikzpicture}
\begin{axis}[ytick = {0}]

\path[name path=up] (axis cs:-1,10) -- (axis cs:1,10);
\path[name path=down] (axis cs:-1,-10) -- (axis cs:1,-10);

\addplot [thick,mark=none] coordinates {(0,0) (1,0)};

\addplot[mark=*,mark size=0.7pt] coordinates {(0,0)};

\addplot [
    color=gray,
    fill=gray, 
    fill opacity=0.2
]
fill between[
    of=up and down,
];

\end{axis}
\end{tikzpicture}
  \caption{\cref{ex:missing_half-line,ex:missing_half-line_infty}}
  \label{fig:missing_half-line}
\end{subfigure}%
\begin{subfigure}{.33\textwidth}
  \centering
  \begin{tikzpicture}
\begin{axis}[ytick = {-1,0,1}]

\path[name path=rightup] (axis cs:0,10) -- (axis cs:1,10);
\path[name path=rightdown] (axis cs:0,-10) -- (axis cs:1,-10);

\addplot [
    name path=leftup,
    thick,
    domain=-1:0, 
    samples=10, 
]
{1};

\addplot [
    name path=leftdown,
    thick,
    domain=-1:0, 
    samples=10, 
]
{-1};

\addplot [thick,mark=none] coordinates {(0,1) (0,10)};
\addplot [thick,mark=none] coordinates {(0,-1) (0,-10)};

\addplot [
    thick,
    color=gray,
    fill=gray, 
    fill opacity=0.2
]
fill between[
    of=leftup and leftdown,
];
\addplot [
    color=gray,
    fill=gray, 
    fill opacity=0.2
]
fill between[
    of=rightup and rightdown,
];

\end{axis}
\end{tikzpicture}
  \caption{\cref{ex:finite_complete}}
  \label{fig:finite_complete}
\end{subfigure}%
\begin{subfigure}{.33\textwidth}
  \centering
  \begin{tikzpicture}
\begin{axis}[ytick = {-1,0,1}]

\path[name path=rightup] (axis cs:-0.1,10) -- (axis cs:1,10);
\path[name path=rightdown] (axis cs:-0.1,-10) -- (axis cs:1,-10);

\addplot [
    name path=leftup,
    thick,
    domain=-1:-0.1, 
    samples=100, 
]
{x^(-1)};

\addplot [
    name path=leftdown,
    thick,
    domain=-1:-0.1, 
    samples=100, 
]
{-x^(-1)};

\addplot [
    color=gray,
    fill=gray, 
    fill opacity=0.2
]
fill between[
    of=leftup and leftdown,
];
\addplot [
    color=gray,
    fill=gray, 
    fill opacity=0.2
]
fill between[
    of=rightup and rightdown,
];

\end{axis}
\end{tikzpicture}
  \caption{\cref{ex:funnel}}
  \label{fig:funnel}
\end{subfigure}%
\caption{Examples of submanifolds $M\subset(-1,1)\times\R$. The shaded area indicates $M$; the thick lines are not part of $M$ (they indicate either the boundary or a removed line). }
\label{fig:regular_examples}
\end{figure}

\begin{example}[A missing half-line]
\label{ex:missing_half-line}
Consider the submanifold $M$ of $(-1,1)\times\R$ obtained by removing a half-line:
\begin{align}
\label{eq:missing_half-line}
M &:= \big((-1,1)\times\R\big) \big\backslash \big([0,1)\times\{0\}\big) 
= \big\{ (x,y)\in(-1,1)\times\R : y\neq0 \text{ whenever } x\geq0 \big\} , 
\end{align}
and consider the flat vertical metric $g_V(x,y) = dy^2$ (i.e.\ $h=1$). 
In this example, $\D$ is \emph{not} regular. 
\end{example}
\begin{proof}
By \cref{lem:vert_regular}, it suffices to show that $ev_x(\Dom\D^*)$ is not equal to $\Dom\D_x^*$ for some $x\in(-1,1)$. 
The critical point is of course $x=0$. 
So consider any element $\eta\in\Dom\D_0^* = H^1((-\infty,0)\cup(0,\infty)) = H^1(-\infty,0) \oplus H^1(0,\infty)$ which is discontinuous at $0$. Clearly $\eta$ extends to an element $\eta\in L^2(\R)$. 
Suppose there exists an element $\psi\in\Dom\D^*$ such that $\psi(0)=\eta$. Since $\psi\in\Gamma_0((-1,1),E_\bullet)$ and $\D^*\psi\in\Gamma_0((-1,1),E_\bullet)$ vary continuously as a function of $x\in(-1,1)$, we must have $\lim_{x\to0^-}\psi(x) = \psi(0)$ and $\lim_{x\to0^-}\partial_y\psi(x) = \partial_y\psi(0)$. 
Since $\psi(x)\in\Dom(-i\partial_y)$ for each $x<0$, this means that $\eta = \psi(0) = \lim_{x\to0^-}\psi(x)$ also lies in $\Dom(-i\partial_y) = H^1(\R) \subset L^2(\R)$. 
However, this cannot be true, since $\eta$ is not continuous on $\R$. 
Thus, by contradiction, it follows that $ev_0(\Dom\D^*)$ does not contain $\eta$, and therefore $\D$ is not regular. 
\end{proof}

\begin{example}[A missing half-line `pushed to infinity']
\label{ex:missing_half-line_infty}
Consider the same manifold $M$ of \cref{eq:missing_half-line}, but replace the flat vertical metric by a `vertically complete' metric $g_V(x,y) = h(x,y) dy^2$ given by the smooth function
\[
h(x,y) := \begin{cases}
          \frac{1}{y^2+e^{1/x}} , & x<0 , \\
          \frac{1}{y^2} , & x\geq0 .
          \end{cases}
\]
Then we know from \cref{lem:complete_regular} that $\D$ is regular. 
\end{example}

\begin{example}[From finite intervals to complete lines]
\label{ex:finite_complete}
Consider the submanifold $M$ of $(-1,1)\times\R$ obtained as the union of $(-1,1)\times(-1,1)$ and $(0,1)\times\R$, equipped with the flat vertical metric $g_V(x,y) = dy^2$. 
We note that for $x\leq0$ the fibre is an incomplete finite interval, while for $x>0$ the fibre is a complete line. 
In this example, $\D$ is \emph{not} regular. 
\end{example}
\begin{proof}
The proof is analogous to the case of \cref{ex:missing_half-line}: one shows that an element $\eta \in \Dom\D_0^* = H^1(-1,1)$ which does not vanish on the boundary cannot be extended to a continuous section in $\Dom\D^*$. 
\end{proof}

Again, if we replace the flat metric in the above example by a complete metric, then $\D$ becomes regular. However, we would like to `fix' the above example without considering a complete metric. Instead, we keep the flat metric but `asymptotically enlarge' the finite intervals as we approach the transition to complete lines. In this way, the transition between incomplete and complete fibres takes place `at infinity', and we obtain a regular operator. 

\begin{example}[The funnel]
\label{ex:funnel}
Consider the open submanifold $M$ of $(-1,1)\times\R$ given by
\[
M := \left\{ (x,y) \in (-1,1)\times\R : |y| x > -1 \right\} ,
\]
equipped with the flat vertical metric $g_V(x,y) = dy^2$. 
We note that for $x<0$ the fibre is an incomplete finite interval, while for $x\geq0$ the fibre is a complete line. 
In this example, $\D$ is regular. 
\end{example}
\begin{proof}
It suffices to show that $\ev_x(\Dom\D^*)$ contains a core for $\D_x^*$. 
We will prove this for the critical point $x=0$, where the transition between complete and incomplete fibres takes place. 
We note that the fibre $\pi^{-1}(0) = \R$ is complete. Thus $\Dom\D_0^* = \Dom\D_0$, and in particular $C_c^\infty(\R)$ is a core for $\D_0^*$. 
Since $C_c^\infty(\R) = \ev_0(C_c^\infty(M))$ is contained in $\ev_0(\Dom\D^*)$, this completes the proof. 
\end{proof}

Finally, merging \cref{ex:missing_half-line_infty,ex:funnel}, we obtain an example which contains a transition from a single incomplete interval to two disjoint complete lines. 

\begin{example}
\label{ex:funnel_missing_half-line}
Consider the open submanifold $M$ of $(-1,1)\times\R$ given by
\[
M := \left\{ (x,y) \in (-1,1)\times\R : |y| x > -1 \; \& \; y\neq0 \text{ whenever } x\geq0\right\} ,
\]
equipped with the vertical metric $g_V(x,y) = h(x,y) dy^2$ given by the smooth function
\[
h(x,y) := \begin{cases}
          \frac{1}{y^2+e^{1/x}} , & x<0 , \\
          \frac{1}{y^2} , & x\geq0 .
          \end{cases}
\]
In this example, $\D$ is again regular (the details are left to the reader). 
\end{example}

\subsubsection{Codimension 2 or higher}
\label{sec:codim_2}

So far, we have seen examples in which the topology of the fibres is allowed to change, but only if this change occurs `at geometric infinity' (see \cref{ex:missing_half-line_infty}). We would like to allow also for such a change of topology without changing the metric, and we will show here that this is possible if the change of topology is obtained for instance by removing a suitable compact subset of codimension 2 or higher. 

Consider the following setup. Let $(N,g_N)$ be a (possibly non-compact) Riemannian manifold, and let $C_N\subset N$ be a compact subset. 
Let $\D_N$ be an elliptic symmetric first-order differential operator on a vector bundle $\bE_N\to N$. 
\begin{lem}
\label{lem:Sobolev_dense}
Suppose that $C_c^\infty(N\backslash C_N)$ is dense in $H_c^1(N)$. 
Then the closure of $\D_N|_{N\backslash C}$ is equal to the closure of $\D_N$. 
\end{lem}
\begin{proof}
Let $\psi \in \Gamma_c^\infty(N,\bE_N)$, and let $\chi\in C_c^\infty(N)$ be such that $\chi = 1$ on a neighbourhood of $\supp\psi$. 
By assumption, there exists a sequence $\chi_n\in C_c^\infty(N\backslash C_N)$ such that $\chi_n\to\chi$ with respect to the $H^1$-norm. 
Hence $\chi_n\psi\to\chi\psi=\psi$ with respect to the graph norm of $\D_N$. Since $\chi_n\psi\in\Gamma_c^\infty(N\backslash C_N,\bE_N) \subset \Dom\D_N|_{N\backslash C}$, this proves the statement. 
\end{proof}

Now let $B$ be another manifold, and consider the obvious submersion $p\colon N\times B \to B$ given by projection onto the second factor. 
Given a closed subset $C_B \subset B$, consider the submanifold $M \subset N\times B$ given by 
\[
M := \big(N\backslash C_N \times B\big) \cup \big(N \times B\backslash C_B\big) .
\]
The projection $p$ then induces a submersion $\pi := p|_M\colon M\to B$. We denote the fibres of this submersion as $M_b := \pi^{-1}(b)$ (for $b\in B$), and we note that we have $M_b = N\backslash C_N$ for $b\in C_B$ and $M_b = N$ for $b\notin C_B$. 

The bundle $\bE_N$ pulls back to a hermitian vector bundle $\bE := \pi^*\bE_N = p^*\bE_N|_{M}$ over $M$. 
We implicitly identify $\Gamma_c^\infty(M,\bE)$ with a subset of $C_c^\infty(B,\Gamma_c^\infty(N,\bE_N))$. 
The operator $\D_N$ induces a vertically elliptic, symmetric, first-order differential operator $\D$ on $\bE$: for $\psi\in\Gamma_c^\infty(M,\bE)$ and $(n,b)\in M$, we have
\[
(\D\psi)(n,b) := \big(\D_N \psi(\cdot,b)\big)(n) ,
\] 
where we implicitly identify $\psi(\cdot,b) \in \Gamma_c^\infty(M_b,\bE|_{M_b})$ with an element in $\Gamma_c^\infty(M_b,\bE_N|_{M_b}) \subset \Gamma_c^\infty(N,\bE_N)$. As usual, we consider the closure of $\D$ as an operator on the Hilbert $C_0(B)$-module $\Gamma_0(B,E_\bullet)$, where $E_b := L^2(M_b,\bE|_{M_b})$. 

\begin{prop}
Suppose that $C_c^\infty(N\backslash C_N)$ is dense in $H_c^1(N)$. 
Then the operator $\D$ on $\Gamma_0(B,E_\bullet)$ is regular. 
\end{prop}
\begin{proof}
By \cref{lem:Sobolev_dense}, we know that $\bar{\D_N|_{N\backslash C_N}} = \bar{\D_N}$ and therefore $(\D_N|_{N\backslash C_N})^* = \D_N^*$. 
Thus $\Dom\D_b^*$ is independent of $b\in B$. 
Consider then any point $b\in B$, and an element $\eta\in\Dom\D_b^*$. 
Picking a function $\phi\in C_c^\infty(B)$ such that $\phi(b)=1$, and defining $\psi(b') := \phi(b') \eta$ for all $b'\in B$, we obtain an element $\psi\in\Dom\D^*$ such that $\psi(b) = \eta$. 
By \cref{lem:vert_regular}, we conclude that $\D$ is regular. 
\end{proof}

\begin{remark}
The assumption that $C_c^\infty(N\backslash C_N)$ is dense in $H_c^1(N)$ is satisfied for instance in the following cases:
\begin{itemize}
\item $C_N$ is a finite union of compact embedded submanifolds of $N$, each of codimension $2$ or higher (see e.g.\ \cite[Lemma 29 \& Proposition 30]{KS17bpre}); 
\item $N\subset\R^n$ is an open subset (equipped with the flat metric), and $C_N$ is a ``compact $d$-set'' of codimension $2$ or higher, i.e.\ a compact subset of Hausdorff-dimension $d \leq n-2$ for which there exist constants $c_1,c_2>0$ such that 
\[
0 < c_1 r^d \leq \mH^d(B_r(x)\cap C_N) \leq c_2 r^d < \infty , 
\]
for all $x\in C_N$ and $0<r<1$, where $\mH^d$ is the $d$-dimensional Hausdorff measure on $\R^n$, and $B_r(x)$ denotes the ball or radius $r$ around $x$ (see \cite[Theorem 2.17]{HM17}). For instance, any $d$-dimensional compact Lipschitz submanifold of $\R^n$ is a compact $d$-set. 
\end{itemize}
\end{remark}

\subsection{Half-closed modules}
\label{sec:half-closed}

In the previous subsection we have considered the regularity property of a vertical differential operator on a submersion. In the following subsection, we will show (assuming regularity) that such an operator defines a half-closed module representing the $\KK$-class of the submersion. 
Before doing so, we need to briefly recall here some definitions and results from \cite{Hil10}. 

Let $A$ be a separable $C^*$-algebra and $B$ a $\sigma$-unital $C^*$-algebra, and let $E$ be a countably generated Hilbert $B$-module.  
For a symmetric operator $\D$ on $E$, we consider the following subspaces of $\End_B(E)$:
\begin{align*}
\Lip(\D) &:= \big\{ T\in\End_B(E) : T\cdot\Dom\D\subset\Dom\D , \text{ and $[\D,T]$ is bounded on $\Dom\D$} \big\} , \\
\Lip^*(\D) &:= \big\{ T\in\Lip(\D) : T\cdot\Dom\D^*\subset\Dom\D \big\} .
\end{align*}
Given a $*$-homomorphism $A\to\End_B(E)$, an operator $T\in\End_B(E)$ is called \emph{locally compact} if $aT$ is compact for every $a\in A$. 

\begin{defn}[{\cite[\S2]{Hil10}}]
\label{defn:Hilsum}
A \emph{half-closed $A$-$B$-module} $(\A,{}_\pi E_B,\D)$ is given by a $\Z_2$-graded countably generated Hilbert $B$-bimodule $E = E^+\oplus E^-$, an odd regular symmetric operator $\D$ on $E$, a $*$-homomorphism $\pi\colon A\to\End_B(E)$, and a dense $*$-subalgebra $\A\subset A$ such that 
\begin{enumerate}
\item $\pi(\A)\subset\Lip^*(\D)$; 
\item $(1+\D^*\D)^{-1}$ is locally compact. 
\end{enumerate}
We will usually suppress the $*$-homomorphism $\pi$ in our notation and simply write $(\A,E_B,\D)$. 
\end{defn}
If furthermore $\D$ is self-adjoint, then $(\A,E_B,\D)$ is called a closed $A$-$B$-module or, more commonly, an \emph{unbounded Kasparov $A$-$B$-module}. Unbounded Kasparov modules were first introduced by Baaj and Julg \cite{BJ83}, who proved that their bounded transforms yield Kasparov modules. This statement was generalised to half-closed modules by Hilsum. 

\begin{thm}[{\cite[Theorem 3.2]{Hil10}}]
\label{thm:Hilsum}
Let $(\A,E_B,\D)$ be a half-closed $A$-$B$-module, and consider a closed extension $\D\subset\hat\D\subset\D^*$. Then the bounded transform $\hat F := F_{\hat\D} = \hat\D(1+\hat\D^*\hat\D)^{-\frac12}$ yields a Kasparov $A$-$B$-module $(A,E_B,\hat F)$, and its class $[\D] := [\hat F] \in \KK(A,B)$ is independent of the choice of the extension $\hat\D$. 
\end{thm}
\begin{remark}
We note that, in particular, the above theorem shows that $[\hat F,a]$ is compact for any $a\in A$, which we will need later on. 
\end{remark}

For any regular operator $\D$, we introduce the notation 
\[
R_\D(\lambda) := (1+\lambda+\D^*\D)^{-1} . 
\]
We recall that we have the integral formula
\begin{align}
\label{eq:integral_formula}
(1+\D^*\D)^{-\frac12} = \frac1\pi \int_0^\infty \lambda^{-1/2} R_\D(\lambda) d\lambda ,
\end{align}
where the integral converges in norm. 

\begin{lem}[{cf.\ \cite[Lemma 7]{Kuc97}}]
\label{lem:integral_formula}
Let $\D$ be a regular operator on $E$. Then for all $\psi\in E$ we have
\[
\frac1\pi \int_0^\infty \lambda^{-1/2} \D R_\D(\lambda) \psi d\lambda = \D(1+\D^*\D)^{-1/2} \psi .
\]
Moreover, for any continuous function $g\colon\R\to\R$ such that $g(x^2)(1+x^2)^{-\frac12}$ is bounded, we also have 
\[
\frac1\pi \int_0^\infty \lambda^{-1/2} g(\D^*\D) R_\D(\lambda) \psi d\lambda = g(\D^*\D)(1+\D^*\D)^{-1/2} \psi .
\]
\end{lem}
\begin{proof}
The first statement is proven for a regular self-adjoint operator $\D$ in \cite[Lemma 7]{Kuc97}. In general, for any regular operator $\D$, we can apply the statement to the self-adjoint operator $\mattwo{0}{\D^*}{\D}{0}$. The proof of the second statement is analogous. 
\end{proof}

The following result shows that two half-closed modules which are `locally bounded' perturbations of each other represent the same class in $\KK$-theory. 

\begin{prop}
\label{prop:loc_bdd}
Consider two half-closed $A$-$B$-modules $(\A,E_B,\D)$ and $(\A,E_B,\D')$ (with the same $*$-homomorphism $\pi\colon A\to\End_B(E)$), and suppose that $\Dom\D\cap\Dom\D'$ is dense in $E$. Assume that for each $a\in\A$, the operator $(\D-\D')a$ extends to a bounded operator on $E$. 
Then $[\D] = [\D'] \in \KK(A,B)$. 
\end{prop}
\begin{proof}
The proof is an adaptation of \cite[Proposition 3.3]{DM19pre} to the setting of half-closed modules. 
Consider the $*$-homomorphism $\til\pi\colon A\to\End_B(E\oplus E)$ given for $t\in[0,1]$ by $\til\pi_t(a) := (a\oplus a) P_t$ in terms of the norm-continuous family of projections 
\[
P_t := \mattwo{\cos^2(\frac{\pi t}2)}{\cos(\frac{\pi t}2)\sin(\frac{\pi t}2)}{\cos(\frac{\pi t}2)\sin(\frac{\pi t}2)}{\sin^2(\frac{\pi t}2)} .
\]
We note that $P_0 = 1\oplus0$ and $P_1 = 0\oplus1$. We will show that we obtain a half-closed $A$-$C([0,1],B)$-module $\big(\A,{}_{\til\pi}C([0,1],E\oplus E)_{C([0,1],B)},\D\oplus\D'\big)$. 
It then follows that its bounded transform is a homotopy between the (bounded) Kasparov modules $(\A, {}_{\pi\oplus0}E\oplus E , F_\D \oplus F_{\D'})$ and $(\A, {}_{0\oplus\pi}E\oplus E , F_\D \oplus F_{\D'})$, which in turn are equal to $(\A, {}_\pi E, F_\D)$ and $(\A, {}_\pi E, F_{\D'})$ (respectively) up to addition of degenerate modules. Thus $[F_\D] = [F_{\D'}]$, which proves the statement. 

The operator $\D\oplus\D'$ is an odd regular symmetric operator on $C([0,1],E\oplus E)$ with locally compact resolvents. To show that $\big(\A,{}_{\til\pi}C([0,1],E\oplus E)_{C([0,1],B)},\D\oplus\D'\big)$ is a half-closed module, we need to show that $\til\pi(\A)\subset\Lip^*(\D\oplus\D')$. 
For any $a\in\A$ we compute 
\[
[ \D \oplus \D' , (a\oplus a) P_t ] = \mattwo{[\D ,a]\cos^2(\frac{\pi t}2)}{(\D a-a\D') \cos(\frac{\pi t}2) \sin(\frac{\pi t}2)}{(\D'a-a\D ) \cos(\frac{\pi t}2)\sin(\frac{\pi t}2)}{[\D',a] \sin^2(\frac{\pi t}2)} .
\]
We observe that $\D a-a\D' = (\D -\D')a + [\D',a]$ is bounded, and similarly for $\D'a-a\D $. 
Hence $[ \D \oplus \D' , (a\oplus a) P_t ]$ is uniformly bounded and norm-continuous in $t$, and we obtain $\til\pi(a) \in \Lip(\D\oplus\D')$. 

It remains to show that $\til\pi(\A)\cdot\Dom(\D\oplus\D')^*\subset\Dom(\D\oplus\D')$. Since $(\A,E_B,\D)$ and $(\A,E_B,\D')$ are half-closed modules, we know that $\A\cdot\Dom\D^*\subset\Dom\D$ and $\A\cdot\Dom(\D')^*\subset\Dom\D'$. Moreover, since $(\D-\D')a$ is bounded for $a\in\A$, we know that $\Dom\D a=\Dom\D'a$, so for any $\psi\in E$ we have $a\psi\in\Dom\D$ if and only if $a\psi\in\Dom\D'$. Thus we have $a\cdot\Dom\D^* \subset \Dom\D\cap\Dom\D'$ and similarly $a\cdot\Dom(\D')^*\subset\Dom\D\cap\Dom\D'$, which proves that 
\begin{align*}
\til\pi(a)\cdot\Dom(\D\oplus\D')^* &\subset P_t \cdot (\Dom\D\cap\Dom\D') \oplus (\Dom\D\cap\Dom\D') \\
&\subset \Dom(\D\oplus\D') . 
\qedhere
\end{align*}
\end{proof}

\subsection{The half-closed module on a submersion}
\label{sec:half-closed_submersion}

We consider as before a submersion $\pi\colon M\to B$, a vertical metric $g_V$ on $M$, and a hermitian vector bundle $\bE\to M$. 

\begin{prop}
\label{prop:vert_cpt_resolvents}
Let $f\in C_c^\infty(M)$, and let $\D$ be (the closure of) a symmetric first-order differential operator on $\bE$, which is vertically elliptic on a neighbourhood of $\supp(f)$. 
Assume that $\D$ is regular as an operator on the Hilbert $C_0(B)$-module $\Gamma_0(B,E_\bullet)$. 
Then the operator $f(1+\D^*\D)^{-\frac12}$ is compact on $\Gamma_0(B,E_\bullet)$. 
\end{prop}
\begin{remark}
The statement also applies to the special case in which $B$ is just a point. In this case, $\D$ is elliptic on $\supp(f)$, and the operator $f(1+\D^*\D)^{-\frac12}$ is compact on the Hilbert space $L^2(M,\bE)$. 
\end{remark}
\begin{proof}
For the compactness of $f (1+\D^*\D)^{-\frac12}$, we need to show that the composition $\Dom\D \xrightarrow{\iota} \Gamma_0(B,E_\bullet) \xrightarrow{f} \Gamma_0(B,E_\bullet)$ is compact, where $\iota$ denotes the domain inclusion. The proof is exactly as in \cite[Propositions 7 \& 11]{KS17bpre} (indeed, though the paper \cite{KS17bpre} focuses on \emph{proper} submersions, the proofs of \cite[Propositions 7 \& 11]{KS17bpre} are local and do not use the properness assumption). 
\end{proof}

\begin{assumption}
\label{ass:vert_op}
Let $\pi\colon M\to B$ be a submersion of manifolds equipped with a vertical metric $g_V$. 
Let $\bE = \bE^+\oplus\bE^- \to M$ be a $\Z_2$-graded hermitian vector bundle, and let $\D$ be an odd, vertically elliptic, symmetric, first-order differential operator on $\bE\to M$. 
We assume that for each $b\in B$, the subspace $\ev_b(\Dom\D^*) \subset \Dom\D_b^*$ is a \emph{core} for $\D_b^*$. 
\end{assumption}

\begin{prop}
\label{prop:vert_half-closed}
Let $\D$ be a vertical operator as in \cref{ass:vert_op}. 
Then $\D$ yields a \emph{half-closed module} $(C_c^\infty(M),\Gamma_0(B,E_\bullet)_{C_0(B)},\D)$ and therefore gives a well-defined class $[\D] := [F_\D] \in \KK_0(C_0(M),C_0(B))$. 
\end{prop}
\begin{proof}
By assumption $\D$ is an odd symmetric operator, and from \cref{lem:vert_regular} we know that $\D$ is regular. 
It is clear that $\Lip(\D)$ contains the smooth compactly supported functions $C_c^\infty(M)$. 
Moreover, since $\D$ is a first-order differential operator, we have for any $f\in C_c^\infty(M)$ that $f\cdot\Dom\D^* \subset \Dom\D$, and therefore $C_c^\infty(M) \subset \Lip^*(\D)$. 
Finally, from \cref{prop:vert_cpt_resolvents} we know that $f(1+\D^*\D)^{-\frac12}$ is compact for any $f\in C_c^\infty(M)$.
Hence $(C_c^\infty(M),\Gamma_0(B,E_\bullet)_{C_0(B)},\D)$ is a half-closed module, and it follows from \cref{thm:Hilsum} that we obtain a class $[\D] := [F_\D] \in \KK_0(C_0(M),C_0(B))$. 
\end{proof}

\subsection{A localised representative}
\label{sec:local_representative}

In this subsection we will describe a \emph{localised representative} $\til F_\D$ for the class $[\D] = [F_\D]$, which is constructed from `localisations' of the unbounded operator $\D$. The construction is due to Higson \cite{Hig89pre} (see also \cite[\S10.8]{Higson-Roe00}), who defined the $\K$-homology class of a first-order symmetric elliptic differential operator $\D$ on an open manifold $M$ which is not necessarily complete (so in particular $\D$ is not necessarily self-adjoint). We will show here that Higson's construction extends to the case of the $\KK$-class of a vertically elliptic operator on a submersion. 

Moreover, it will be useful later on to allow for some additional flexibility in Higson's construction, by rescaling each localisation of $\D$ independently. To show that this rescaling is allowed, we first prove the following lemma. 

\begin{lem}
\label{lem:bdd_transf_rescaled}
Let $B$ be a $C^*$-algebra, and let $\D$ be a regular self-adjoint operator on a Hilbert $B$-module $E$. Let $a\in\End_B(E)$ such that $a(\D\pm i)^{-1}$ is compact. Then for any $\alpha>0$, the operator $a (F_\D - F_{\alpha\D})$ is compact. 
\end{lem}
\begin{proof}
Since $a(\D\pm i)^{-1}$ is compact, and since the functions $x\mapsto(x\pm i)^{-1}$ generate $C_0(\R)$, we see that $ag(D)$ is compact for any $g\in C_0(\R)$. The statement then follows because the function $g_\alpha\colon x \mapsto x(1+x^2)^{-\frac12} - \alpha x(1+\alpha^2x^2)^{-\frac12}$ lies in $C_0(\R)$. 
\end{proof}

Throughout the remainder of this section, we consider the following setting. 

\begin{assumption}
\label{ass:Higson}
Let $\D$ be as in \cref{ass:vert_op}, and let $\{U_j\}_{j\in\N}$ be a locally finite cover of open precompact subsets of $M$. 
We consider compactly supported, smooth, real-valued functions $\{\chi_j\}_{j\in\N}$ and $\{\phi_j\}_{j\in\N}$ satisfying the following properties:
\begin{enumerate}
\item $\{\chi_j^2\}$ is a partition of unity for $\{U_j\}$; 
\item $\phi_j|_{U_j}=1$ for all $j\in\N$.
\end{enumerate}
\end{assumption}

For each $j$, consider the first-order differential operator $\D_j := \phi_j \D \phi_j$. Since $\D_j$ is a compactly supported, symmetric, first-order differential operator, we observe that $\D_j$ is in fact \emph{essentially self-adjoint}. We rescale each $\D_j$ by a positive number $\alpha_j$, and consider the bounded transforms $F_{\alpha_j\D_j} = \alpha_j\D_j (1+\alpha_j^2\D_j^2)^{-\frac12}$. 

\begin{defn}
\label{defn:Higson}
Given a sequence $\{\alpha_j\}_{j\in\N} \subset (0,\infty)$ of positive numbers, 
we define the \emph{localised representative} of $\D$ as 
\begin{align*}
\til F_\D(\alpha) := \sum_j \chi_j F_{\alpha_j\D_j} \chi_j .
\end{align*}
\end{defn}

The inequality $\sum_{j=0}^k \chi_j F_{\alpha_j\D_j} \chi_j \leq 1$ shows that the partial sums  are uniformly bounded. Moreover, since the partition of unity is locally finite, we have for any compactly supported $\psi \in \Gamma_0(B,E_\bullet)$ that $\til F_\D \psi$ is a finite (hence convergent) series. Thus $\til F_\D(\alpha)$ is well-defined on $\Gamma_0(B,E_\bullet)$ as a strongly convergent series. 

We will show that $\til F_\D(\alpha)$ defines a Kasparov module, and that its class is equal to the class of the bounded transform $F_\D$. It suffices to show that the difference $\til F_\D(\alpha) - F_\D$ is locally compact. It then follows in particular that the Kasparov class $[\til F_\D(\alpha)]$ is independent of the choices made in the construction. 

\begin{lem}
\label{lem:R_S_cpt_order}
Consider a smooth self-adjoint endomorphism $T\in\Gamma_c^\infty(M,\End\bE) \cap \Lip(\D)$, and suppose that $\D$ is (the closure of) a symmetric first-order differential operator on $\bE$, which is vertically elliptic on a neighbourhood of $\supp(T)$. 
If $\D$ is regular on $\Gamma_0(B,E_\bullet)$, then the operators $[\D R_\D(\lambda),T]$ and $(1+\lambda)^{\frac12} [R_\D(\lambda),T]$ are compact and of order $\mO(\lambda^{-1})$. 
\end{lem}
\begin{proof}
We have
\begin{align*}
[R_\D(\lambda),T] &= \big[ (1+\lambda+\D^*\D)^{-1} , T \big]  = - (1+\lambda+\D^*\D)^{-1} [\D^*\D,T] (1+\lambda+\D^*\D)^{-1}  .
\end{align*}
Rewriting $[\D^*\D,T] = [\D^*,T] \D + \D^* [\D,T]$, we obtain 
\[
[R_\D(\lambda),T]  = - R_\D(\lambda)^{\frac12} \Big( R_\D(\lambda)^{\frac12} [\D^*,T] \big( \D R_\D(\lambda)^{\frac12} \big) + \big( R_\D(\lambda)^{\frac12} \D^* \big) [\D,T] R_\D(\lambda)^{\frac12} \Big) R_\D(\lambda)^{\frac12}  .
\]
We note that $[\D^*,T] = [\D,T]$ is a smooth endomorphism with $\supp([\D,T]) \subset \supp(T)$. 
Since $\D$ is vertically elliptic on a neighbourhood of $\supp(T)$, we know from \cref{prop:vert_cpt_resolvents} that $[\D,T] R_\D(\lambda)^{\frac12}$ and $R_\D(\lambda)^{\frac12} [\D^*,T]$ are compact. Hence $[R_\D(\lambda),T]$ is compact and of order $\mO(\lambda^{-\frac32})$. Moreover, we see that $\D [R_\D(\lambda),T]$ is also well-defined, compact, and of order $\mO(\lambda^{-1})$. Thus we conclude that 
\[
[\D R_\D(\lambda),T]  = [\D,T] R_\D(\lambda)  + \D [R_\D(\lambda),T]  
\]
is compact and of order $\mO(\lambda^{-1})$. 
\end{proof}

\begin{lem}
\label{lem:local_comparison}
Consider the setting of \cref{ass:Higson}. 
Let $f\in C_c^\infty(M)$ and $\phi\in C_c^\infty(M,\R)$ be such that $\phi=1$ on a neighbourhood of $\supp f$, 
so that $\D_\phi := \phi\D\phi$ is an essentially self-adjoint operator which agrees with $\D$ on $\supp f$. 
Then $f(F_\D-F_{\D_\phi})$ is compact on $\Gamma_0(B,E_\bullet)$. 
\end{lem}
\begin{proof}
The proof is very similar to the argument of \cite[Lemma 3.1]{Hil10}. 
Since $f(F_\D-F_\D^*)$ is compact (by \cref{prop:vert_half-closed,thm:Hilsum}), it suffices to show that $f(F_\D^*-F_{\D_\phi})$ is compact. 
We can rewrite 
\begin{align*}
f(F_\D^*-F_{\D_\phi}) 
&= f \left( \bar{(1+\D^*\D)^{-\frac12} \D^*} - \bar{\D_\phi (1+\D_\phi^2)^{-\frac12}} \right) \\
&= f \left( \bar{\D^* (1+\D\D^*)^{-\frac12}} - \bar{\D_\phi (1+\D_\phi^2)^{-\frac12}} \right) . 
\end{align*}
Using \cref{lem:integral_formula}, we have for any $\psi\in \Gamma_0(B,E_\bullet)$ that 
\[
f(F_\D^*-F_{\D_\phi}) \psi = \frac1\pi \int_0^\infty \lambda^{-\frac12} T(\lambda) \psi d\lambda ,
\]
where 
\[
T(\lambda) := f \left( \D^* (1+\lambda+\D\D^*)^{-1} - (1+\lambda+\D_\phi^2)^{-1} \D_\phi \right) .
\]
We claim that $T(\lambda)$ is a compact operator on $\Gamma_0(B,E_\bullet)$, and that $\|T(\lambda)\| = \mO(\lambda^{-1})$ as $\lambda\to\infty$. It then follows that $\frac1\pi \int_0^\infty \lambda^{-\frac12} T(\lambda) d\lambda$ is in fact a \emph{norm}-convergent integral of compact operators, which proves the statement. 

To prove the claim, we rewrite 
\begin{align*}
T(\lambda) 
&= f \bar{(1+\lambda+\D_\phi^2)^{-1} (1+\lambda+\D_\phi^2)} \; \bar{\D^* (1+\lambda+\D\D^*)^{-1}} \\
&\qquad- f \bar{(1+\lambda+\D_\phi^2)^{-1} \D_\phi} \; \bar{(1+\lambda+\D\D^*) (1+\lambda+\D\D^*)^{-1}} \\
&= f \bar{(1+\lambda+\D_\phi^2)^{-1} \D_\phi} \; \bar{(\D_\phi - \D) \D^* (1+\lambda+\D\D^*)^{-1}} \\
&\qquad+ (1+\lambda) f (1+\lambda+\D_\phi^2)^{-1} \bar{(\D^* - \D_\phi) (1+\lambda+\D\D^*)^{-1}} .
\end{align*}
We note that the operators on the last line are still well-defined. For instance, we have $\Ran\big(\D^* (1+\lambda+\D\D^*)^{-1}\big) \subset \Dom\D \subset \Dom\D_\phi$, so that $\bar{(\D_\phi - \D) \D^* (1+\lambda+\D\D^*)^{-1}}$ is a well-defined bounded operator. We also note that $\phi\cdot\Dom\D^*\subset\Dom\D$, so that $\D_\phi$ is well-defined on $\Dom\D^*$. 

Since $\phi|_{\supp f} = 1$, we note that $f(\D_\phi-\D) = f(\D^*-\D_\phi) = 0$. Hence we find that 
\begin{multline*}
T(\lambda) = \big[ f , \bar{(1+\lambda+\D_\phi^2)^{-1} \D_\phi} \big] \bar{(\D_\phi - \D) \D^* (1+\lambda+\D\D^*)^{-1}} \\
+ (1+\lambda)  \big[ f , (1+\lambda+\D_\phi^2)^{-1} \big] \bar{(\D^* - \D_\phi) (1+\lambda+\D\D^*)^{-1}} .
\end{multline*}
For the first term we note that $\bar{(\D_\phi - \D) \D^* (1+\lambda+\D\D^*)^{-1}}$ is bounded (and of order $\mO(\lambda^0)$), and from \cref{lem:R_S_cpt_order} we know that $ \big[ f , \bar{(1+\lambda+\D_\phi^2)^{-1} \D_\phi} \big]$ is compact and of order $\mO(\lambda^{-1})$. 
For the second term, we find similarly that $\bar{(\D^* - \D_\phi) (1+\lambda+\D\D^*)^{-1}} = \bar{(\D^* - \D_\phi) (1+\lambda+\D\D^*)^{-\frac12}} \; (1+\lambda+\D\D^*)^{-\frac12}$ is bounded and of order $\mO(\lambda^{-\frac12})$, and that $ \big[ f , (1+\lambda+\D_\phi^2)^{-1} \big]$ is compact and of order $\mO(\lambda^{-\frac32})$. 
Hence we see that $T(\lambda)$ is compact and of order $\mO(\lambda^{-1})$, as claimed. 
\end{proof}

\begin{thm}
\label{thm:local_rep_KK}
Consider the setting of \cref{ass:Higson}, and let $\til F_\D(\alpha)$ be the localised representative of $\D$ (as constructed in \cref{defn:Higson}). 
Then for any $f\in C_c(M)$, the operator $f(\til F_\D(\alpha)-F_\D)$ is compact. 
Hence $[\til F_\D(\alpha)] = [F_{\D}] \in \KK(C_0(M),C_0(B))$. 
In particular, the class $[\til F_\D(\alpha)]$ is independent of the choices made in the construction. 
\end{thm}
\begin{proof}
Since $\D_j$ is vertically elliptic on (a neighbourhood of) the support of $\chi_j$, we know from \cref{prop:vert_cpt_resolvents} that $\chi_j (\D_j\pm i)^{-1}$ is compact. 
Hence we can apply \cref{lem:bdd_transf_rescaled}, and we see that $\chi_j(F_{\alpha_j\D_j}-F_{\D_j})$ is compact. 
By \cref{thm:Hilsum}, the commutator $[F_\D,\chi_j]$ is compact. 
Furthermore, from \cref{lem:local_comparison} we know that $\chi_j(F_{\D_j}-F_{\D})$ is compact. 
Since the partition of unity is locally finite, we know that $f\til F_\D(\alpha)$ is given by a finite sum, and therefore
\begin{align*}
f(\til F_\D(\alpha)-F_{\D}) 
&= \sum_j f \chi_j F_{\alpha_j\D_j} \chi_j - f F_{\D} 
\sim \sum_j f \big( \chi_j F_{\D_j} \chi_j - \chi_j^2 F_{\D} \big) \\
&\sim \sum_j f \chi_j ( F_{\D_j} - F_{\D} ) \chi_j 
\sim 0 .
\qedhere
\end{align*}
\end{proof}

\section{Local positivity}
\label{sec:local_positivity}

This section contains the technical part of the proof of our main result, \cref{thm:submersion_Kasp_prod}. 
The goal in this section is to show that a `local positivity condition' for two first-order differential operators implies a `local positivity condition' for their bounded transforms (for the precise statement, see \cref{prop:local_positivity} below). 
We consider the following setting. 

\begin{assumption}
\label{ass:local_positivity}
Let $\rho\in C_c^\infty(M,\R)$ such that $\|\rho\|\leq1$, and let $\D$ and $S$ be two odd essentially self-adjoint first-order differential operators on a $\Z_2$-graded hermitian vector bundle $\bE\to M$. We view (the closures of) $\D$ and $S$ as self-adjoint operators on the Hilbert space $L^2(M,\bE)$, and we make the following assumptions:
\begin{enumerate}
\item $\Dom\D \subset \Dom S$;
\item $\D$ is elliptic on a neighbourhood of $\supp \rho$. 
\end{enumerate}
\end{assumption}
We will apply the results of this section to the case where we have a submersion $M\to B$ and where $S$ is a vertical operator, so the reader may keep this case in mind. 

By the closed graph theorem, assumption (1) implies that $S (1+\D^2)^{-\frac12}$ is bounded. 
Furthermore, since $\D$ has smooth coefficients, $[\D,\rho]$ preserves the domain of $\D$ and hence also $S [\D,\rho] (1+\D^2)^{-\frac12}$ is bounded. 
We will use these facts throughout this section. 
For $\lambda,\mu\in[0,\infty)$, we define 
\begin{align*}
R_\D(\lambda) &:= (1+\lambda+\D^2)^{-1} , & 
R_S(\mu) &:= (1+\mu+S^2)^{-1} . 
\end{align*}
Thus for any $\psi\in L^2(M,\bE)$ we have by \cref{lem:integral_formula} that $F_\D \psi = \frac1\pi \int_0^\infty \lambda^{-\frac12} \D R_\D(\lambda) \psi d\lambda$, and similarly for $F_S$. 
We introduce the following bounded operators: 
\begin{align*}
K(\lambda,\mu) &:= (1+\mu) R_S(\mu) [\D R_\D(\lambda),\rho^2]^* S R_S(\mu) , & 
B_1(\lambda) &:= \D R_\D(\lambda) \rho [S,\rho] , \\
B_2(\lambda) &:= (1+\lambda) R_\D(\lambda) [\D,\rho]^* S \rho R_\D(\lambda) , & 
B_3(\lambda) &:= \D R_\D(\lambda) \rho S [\D,\rho] \D R_\D(\lambda) , \\
M_1(\lambda,\mu) &:= \D R_\D(\lambda) S R_S(\mu) , & 
M_2(\lambda,\mu) &:= \D R_\D(\lambda) \sqrt{1+\mu} R_S(\mu) , \\
M_3(\lambda,\mu) &:= \sqrt{1+\lambda} R_\D(\lambda) S R_S(\mu) , & 
M_4(\lambda,\mu) &:= \sqrt{1+\lambda} R_\D(\lambda) \sqrt{1+\mu} R_S(\mu) .
\end{align*}
Moreover, we consider the quadratic form $Q$ defined for $\psi \in \Dom\D$ by
\[
Q(\psi) := 2 \Re \la \D \psi | S \psi \ra . 
\]

In this section, we will study the positivity of the operator $\chi [F_\D,F_S]_\pm \chi$ for some $\chi \in C_c^\infty(M,\R)$. 
Applying \cref{lem:integral_formula} twice, we can rewrite 
\begin{align}
\label{eq:comm_integral}
\big\la \psi \bigmvert \chi [F_\D , F_S]_\pm \chi \psi \big\ra 
&= 2 \Re \big\la \chi \psi \bigmvert F_\D F_S \chi \psi \big\ra \nonumber\\
&= \frac{1}{\pi^2} \int_0^\infty \int_0^\infty \lambda^{-\frac12} \mu^{-\frac12} 2 \Re \big\la \chi \psi \bigmvert \D R_\D(\lambda) S R_S(\mu) \chi \psi \big\ra d\mu d\lambda .
\end{align}
Our first task is to study the integrand on the right-hand-side. Via a straightforward but somewhat tedious calculation, we will rewrite this integrand in terms of the operators $K(\lambda,\mu)$, $B_l(\lambda)$, and $M_m(\lambda,\mu)$ defined above. 

\begin{lem}
\label{lem:positivity_local_expr}
For any $\psi\in L^2(M,\bE)$ we have 
\begin{multline*}
2 \Re \big\la \rho^2 \psi \bigmvert \D R_\D(\lambda) S R_S(\mu) \psi \big\ra 
= 2 \Re \big\la \psi \bigmvert K(\lambda,\mu) \psi \big\ra + \sum_{m=1}^4 Q\big(\rho M_m(\lambda,\mu)\psi\big) \\
- 2 \sum_{l=1}^3 \Re \big\la B_l(\lambda) \sqrt{1+\mu} R_S(\mu) \psi \bigmvert \sqrt{1+\mu} R_S(\mu) \psi \big\ra 
- 2 \sum_{l=1}^3 \Re \big\la B_l(\lambda) S R_S(\mu) \psi \bigmvert S R_S(\mu) \psi \big\ra .
\end{multline*}
\end{lem}
\begin{proof}
First, we calculate
\begin{align}
\label{eq:local_pos_1}
\big\la \rho^2 \psi \bigmvert \D R_\D(\lambda) S R_S(\mu) \psi \big\ra 
&= \big\la \rho^2 (1+\mu+S^2) R_S(\mu) \psi \bigmvert \D R_\D(\lambda) S R_S(\mu) \psi \big\ra \nonumber\\
&= (1+\mu) \big\la \rho^2 R_S(\mu) \psi \bigmvert \D R_\D(\lambda) S R_S(\mu) \psi \big\ra \nonumber\\
&\quad+ \big\la \rho^2 S^2 R_S(\mu) \psi \bigmvert \D R_\D(\lambda) S R_S(\mu) \psi \big\ra \nonumber\\
&= (1+\mu) \big\la \D R_\D(\lambda) \rho^2 R_S(\mu) \psi \bigmvert S R_S(\mu) \psi \big\ra \nonumber\\
&\quad+ \big\la \rho^2 S^2 R_S(\mu) \psi \bigmvert \D R_\D(\lambda) S R_S(\mu) \psi \big\ra \nonumber\\
&= (1+\mu) \big\la [\D R_\D(\lambda),\rho^2] R_S(\mu) \psi \bigmvert S R_S(\mu) \psi \big\ra \nonumber\\
&\quad+ (1+\mu) \big\la \rho \D R_\D(\lambda) R_S(\mu) \psi \bigmvert \rho S R_S(\mu) \psi \big\ra \nonumber\\
&\quad+ \big\la \rho S^2 R_S(\mu) \psi \bigmvert \rho \D R_\D(\lambda) S R_S(\mu) \psi \big\ra \nonumber\\
&= \big\la \psi \bigmvert K(\lambda,\mu) \psi \big\ra + (1+\mu) \big\la \rho \D R_\D(\lambda) R_S(\mu) \psi \bigmvert \rho S R_S(\mu) \psi \big\ra \nonumber\\
&\quad+ \big\la \rho S^2 R_S(\mu) \psi \bigmvert \rho \D R_\D(\lambda) S R_S(\mu) \psi \big\ra .
\end{align}
Considering $\xi = \sqrt{1+\mu} R_S(\mu) \psi$ or $\xi = S R_S(\mu) \psi$, we can rewrite 
\begin{align*}
\big\la \rho S \xi \bigmvert \rho \D R_\D(\lambda) \xi \big\ra 
&= - \big\la [S,\rho] \xi \bigmvert \rho \D R_\D(\lambda) \xi \big\ra 
+ \big\la \rho \xi \bigmvert S \rho \D R_\D(\lambda) \xi \big\ra \\
&= - \big\la [S,\rho] \xi \bigmvert \rho \D R_\D(\lambda) \xi \big\ra 
+ \big\la \rho (1+\lambda+\D^2) R_\D(\lambda) \xi \bigmvert S \rho \D R_\D(\lambda) \xi \big\ra \\
&= - \big\la [S,\rho] \xi \bigmvert \rho \D R_\D(\lambda) \xi \big\ra 
+ (1+\lambda) \big\la \rho R_\D(\lambda) \xi \bigmvert S \rho \D R_\D(\lambda) \xi \big\ra \\
&\quad+ \big\la \rho \D^2 R_\D(\lambda) \xi \bigmvert S \rho \D R_\D(\lambda) \xi \big\ra \\
&= - \big\la [S,\rho] \xi \bigmvert \rho \D R_\D(\lambda) \xi \big\ra 
- (1+\lambda) \big\la S \rho R_\D(\lambda) \xi \bigmvert [\D,\rho] R_\D(\lambda) \xi \big\ra \\
&\quad+ (1+\lambda) \big\la S \rho R_\D(\lambda) \xi \bigmvert \D \rho R_\D(\lambda) \xi \big\ra 
- \big\la [\D,\rho] \D R_\D(\lambda) \xi \bigmvert S \rho \D R_\D(\lambda) \xi \big\ra \\
&\quad+ \big\la \D \rho \D R_\D(\lambda) \xi \bigmvert S \rho \D R_\D(\lambda) \xi \big\ra \\
&= - \big\la B_1(\lambda) \xi \bigmvert \xi \big\ra 
- \big\la B_2(\lambda) \xi \bigmvert \xi \big\ra \\
&\quad+ (1+\lambda) \big\la S \rho R_\D(\lambda) \xi \bigmvert \D \rho R_\D(\lambda) \xi \big\ra 
- \big\la B_3(\lambda) \xi \bigmvert \xi \big\ra \\
&\quad+ \big\la \D \rho \D R_\D(\lambda) \xi \bigmvert S \rho \D R_\D(\lambda) \xi \big\ra .
\end{align*}
Inserting the definition of $Q$, we find
\begin{align*}
2 \Re \big\la \rho S \xi \bigmvert \rho \D R_\D(\lambda) \xi \big\ra 
&= - 2 \sum_{l=1}^3 \Re \big\la B_l(\lambda) \xi \bigmvert \xi \big\ra 
+ (1+\lambda) 2 \Re \big\la S \rho R_\D(\lambda) \xi \bigmvert \D \rho R_\D(\lambda) \xi \big\ra \\
&\quad+ 2 \Re \big\la \D \rho \D R_\D(\lambda) \xi \bigmvert S \rho \D R_\D(\lambda) \xi \big\ra \\
&= - 2 \sum_{l=1}^3 \Re \big\la B_l(\lambda) \xi \bigmvert \xi \big\ra + Q\big(\sqrt{1+\lambda} \rho R_\D(\lambda) \xi\big) + Q\big( \rho \D R_\D(\lambda) \xi\big) . 
\end{align*}
Inserting this expression into \cref{eq:local_pos_1}, we find that 
\begin{align*}
2 \Re &\big\la \rho^2 \psi \bigmvert \D R_\D(\lambda) S R_S(\mu) \psi \big\ra \\
&= 2 \Re \big\la \psi \bigmvert K(\lambda,\mu) \psi \big\ra 
+ (1+\mu) 2 \Re \big\la \rho \D R_\D(\lambda) R_S(\mu) \psi \bigmvert \rho S R_S(\mu) \psi \big\ra \\
&\quad+ 2 \Re \big\la \rho S^2 R_S(\mu) \psi \bigmvert \rho \D R_\D(\lambda) S R_S(\mu) \psi \big\ra \\
&= 2 \Re \big\la \psi \bigmvert K(\lambda,\mu) \psi \big\ra - 2 \sum_{l=1}^3 \Re \big\la B_l(\lambda) \sqrt{1+\mu} R_S(\mu) \psi \bigmvert \sqrt{1+\mu} R_S(\mu) \psi \big\ra \\
&\quad+ Q\big(\sqrt{1+\lambda} \rho R_\D(\lambda) \sqrt{1+\mu} R_S(\mu) \psi\big) + Q\big( \rho \D R_\D(\lambda) \sqrt{1+\mu} R_S(\mu) \psi\big) \\
&\quad- 2 \sum_{l=1}^3 \Re \big\la B_l(\lambda) S R_S(\mu) \psi \bigmvert S R_S(\mu) \psi \big\ra \\
&\quad+ Q\big(\sqrt{1+\lambda} \rho R_\D(\lambda) S R_S(\mu) \psi\big) + Q\big( \rho \D R_\D(\lambda) S R_S(\mu) \psi\big) .
\qedhere
\end{align*}
\end{proof}

Our aim is to control the integrals of each of the terms in the result of \cref{lem:positivity_local_expr}. For the first term, we will show that it gives rise to a compact operator. For the other terms, we will show that we can obtain suitable lower bounds. We will make frequent use of the following lemma. 

\begin{lem}[{\cite[Proposition A.1]{Les05}}]
\label{lem:interpolation}
Let $\mH$ be a separable Hilbert space, let $P$ be an invertible positive self-adjoint operator on $\mH$, and let $T$ be a symmetric operator on $\mH$ with $\Dom P \subset \Dom T$. 
If $TP^{-1}$ is bounded, then the densely defined operator $P^{-\frac12}TP^{-\frac12}$ extends to a bounded operator on $\mH$, and $\|P^{-\frac12}TP^{-\frac12}\| \leq \|TP^{-1}\|$. 
\end{lem}

\begin{lem}
\label{lem:Lipschitz}
Let $\D$ be a self-adjoint operator on a Hilbert space $\mH$, and let $a=a^*\in\Lip(\D)$. 
\begin{enumerate}
\item If $[\D,a]$ preserves $\Dom\D$, then $[(1+\D^2)^{\frac12},a]$ is bounded. 
\item The commutator $[(1+\D^2)^{\frac12},a]$ is bounded if and only if $\D[F_\D,a] $ is bounded. 
\end{enumerate}
\end{lem}
\begin{proof}
\begin{enumerate}
\item 
Rewriting and applying \cref{eq:integral_formula}, we obtain
\begin{align*}
\big[ (1+\D^2)^{\frac12} ,a \big] 
&= - (1+\D^2)^{\frac12} \big[ (1+\D^2)^{-\frac12} , a \big] (1+\D^2)^{\frac12} \\
&= - \frac1\pi \int_0^\infty \lambda^{-1/2} (1+\D^2)^{\frac12} \big[ R_\D(\lambda) , a \big] (1+\D^2)^{\frac12} d\lambda \\
&= \frac1\pi \int_0^\infty \lambda^{-1/2} (1+\D^2)^{\frac12} R_\D(\lambda) \big[ \D^2 , a \big] R_\D(\lambda) (1+\D^2)^{\frac12} d\lambda 
\end{align*}
Since $[\D,a]$ preserves $\Dom\D$, we see that $[\D^2,a] (1+|\D|)^{-1} = \big(\D[\D,a]+[\D,a]\D\big) (1+|\D|)^{-1}$ is well-defined and bounded. Since $i[\D^2,a]$ is symmetric, we know from \cref{lem:interpolation} that also $(1+|\D|)^{-\frac12} i[\D^2,a] (1+|\D|)^{-\frac12}$ is bounded, and that we have 
\[
\big\| (1+|\D|)^{-\frac12} i[\D^2,a] (1+|\D|)^{-\frac12} \big\| \leq \big\| [\D^2,a] (1+|\D|)^{-1} \big\| . 
\]
Hence we obtain the operator inequalities
\begin{align*}
\pm &(1+\D^2)^{\frac12} R_\D(\lambda) i \big[ \D^2 , a \big] R_\D(\lambda) (1+\D^2)^{\frac12} \\
&\leq \big\| [\D^2,a] (1+|\D|)^{-1} \big\| \, (1+\D^2)^{\frac12} R_\D(\lambda) (1+|\D|) R_\D(\lambda) (1+\D^2)^{\frac12} \\
&= \big\| [\D^2,a] (1+|\D|)^{-1} \big\| \, (1+\D^2) (1+|\D|) R_\D(\lambda)^2 \\
&\leq \big\| [\D^2,a] (1+|\D|)^{-1} \big\| \, (1+|\D|) R_\D(\lambda) .
\end{align*}
Inserting this into the integral expresssion for $\big[ (1+\D^2)^{\frac12} ,a \big]$, we find
\begin{align*}
\pm i \big[ (1+\D^2)^{\frac12} ,a \big] 
&\leq \frac1\pi \big\| [\D^2,a] (1+|\D|)^{-1} \big\| \int_0^\infty \lambda^{-1/2} (1+|\D|) R_\D(\lambda) d\lambda . 
\end{align*}
More precisely, since by \cref{lem:integral_formula} the integral converges only strongly (and not in norm), we have for any $\psi\in\Dom\D$ that 
\begin{align*}
\pm \big\la \psi \bigmvert i \big[ (1+\D^2)^{\frac12} ,a \big] \psi \big\ra 
&\leq \frac1\pi \big\| [\D^2,a] (1+|\D|)^{-1} \big\| \int_0^\infty \lambda^{-1/2} \big\la \psi \bigmvert (1+|\D|) R_\D(\lambda) \psi \big\ra d\lambda \\
&= \big\| [\D^2,a] (1+|\D|)^{-1} \big\| \, \big\la \psi \bigmvert (1+|\D|) (1+\D^2)^{-\frac12} \psi \big\ra . 
\end{align*}
Since $(1+|\D|) (1+\D^2)^{-\frac12}$ is bounded, we conclude that $\big[ (1+\D^2)^{\frac12} ,a \big]$ is densely defined and bounded, and therefore it extends to a bounded operator on all of $\mH$. 

\item 
First, we rewrite 
\begin{align*}
\D F_\D 
&= \D^2 (1+\D^2)^{-\frac12} 
= \D^2 (1+\D^2)^{-1} (1+\D^2)^{\frac12} \\
&= \big( 1 - (1+\D^2)^{-1} \big) (1+\D^2)^{\frac12} 
= (1+\D^2)^{\frac12} - (1+\D^2)^{-\frac12} .
\end{align*}
Hence we see that 
\[
\D[F_\D,a] = [\D F_\D,a] - [\D,a] F_\D = \big[ (1+\D^2)^{\frac12} ,a \big] - \big[ (1+\D^2)^{-\frac12} ,a \big] - [\D,a] F_\D .
\]
Since the last two terms on the right-hand-side are bounded, the second statement follows. 
\qedhere
\end{enumerate}
\end{proof}

\begin{lem}
\label{lem:integral_K}
The integral $K := \int_0^\infty \int_0^\infty \mu^{-\frac12} \lambda^{-\frac12} K(\lambda,\mu) d\lambda d\mu$ defines a compact operator $K$. 
\end{lem}
\begin{proof}
We know from \cref{lem:R_S_cpt_order} that 
$[\D R_\D(\lambda),\rho^2]$ is compact and of order $\mO(\lambda^{-1})$. 
Therefore the integral $\int_0^\infty \lambda^{-\frac12} [\D R_\D(\lambda),\rho^2] d\lambda$ converges in norm to the compact operator $[F_\D,\rho^2]$. Hence we obtain 
\[
K = \int_0^\infty \mu^{-\frac12} (1+\mu) R_S(\mu) [F_\D,\rho^2]^* S R_S(\mu) d\mu 
\]
Since $[\D,\rho^2]$ preserves $\Dom\D$, we know from \cref{lem:Lipschitz} that $\D[F_\D,\rho^2]$ is bounded. By \cref{ass:local_positivity} we have $\Dom\D\subset\Dom S$, and therefore $S [F_\D,\rho^2]$ is bounded as well. 
Consequently, also $[F_\D,\rho^2]^* S = \big( S [F_\D,\rho^2] \big)^*$ is bounded. Hence the integrand on the right-hand-side of the above expression is of order $\mO(\mu^{-\frac32})$, so the integral converges in norm. Since the integrand is compact, this proves that $K$ is compact. 
\end{proof}

\begin{lem}
\label{lem:lower_bound_integral}
For $\psi \in L^2(M,\bE)$, consider the integral $I(\psi)$ given by 
\[
% I(\psi) := 
2 \Re \Big\la \psi \Bigmvert \sum_{l=1}^3 \int_0^\infty \int_0^\infty (\lambda\mu)^{-\frac12} \Big( (1+\mu) R_S(\mu) B_l(\lambda) R_S(\mu) + S R_S(\mu) B_l(\lambda) S R_S(\mu) \Big) \psi d\lambda d\mu \Big\ra .
\]
Then there exists a constant $C = C(\D,S,\rho) > 0$ such that for all $\psi \in L^2(M,\bE)$ we have 
\[
\pm I(\psi) \geq - C \la\psi|\psi\ra .
\]
Moreover, if we replace $S$ by $\alpha S$ for some $\alpha>0$, then $C$ is replaced by $\alpha C$. 
\end{lem}
\begin{proof}
We start by deriving some norm estimates for the operators $B_l(\lambda)$. 
First, for $l=1$, we observe that $R_\D(\lambda)^{-\frac12} B_1(\lambda)$ and $R_\D(\lambda)^{-\frac12} B_1(\lambda)^*$ are bounded operators (where for the latter we use that $[S,\rho]\rho$ is smooth and therefore preserves the domain of $\D$). 
Since also $R_\D(\lambda)^{-\frac12} (1+|\D|)^{-1}$ is bounded, we obtain a bounded operator
\[
R_\D(\lambda)^{-\frac12} \left( B_1(\lambda) + B_1(\lambda)^* \right) R_\D(\lambda)^{-\frac12} (1+|\D|)^{-1} .
\]
Applying \cref{lem:interpolation} to the positive invertible operator $P = 1+|\D|$ and the symmetric operator $T = R_\D(\lambda)^{-\frac12} \left( B_1(\lambda) + B_1(\lambda)^* \right) R_\D(\lambda)^{-\frac12}$, we find that also 
\[
\til B_1(\lambda) := (1+|\D|)^{-\frac12} R_\D(\lambda)^{-\frac12} \left( B_1(\lambda) + B_1(\lambda)^* \right) R_\D(\lambda)^{-\frac12} (1+|\D|)^{-\frac12} 
\]
is bounded. We note that we can write
\begin{align*}
B_1(\lambda) + B_1(\lambda)^* &= (1+|\D|)^{\frac12} R_\D(\lambda)^{\frac12} \til B_1(\lambda) (1+|\D|)^{\frac12} R_\D(\lambda)^{\frac12} .
\end{align*}
Since for any self-adjoint endomorphisms $B$ and $S$ we have $\pm B \leq \|B\|\cdot\Id$ and therefore $\pm SBS \leq \|B\| S^2$, we obtain:
\begin{align*}
\pm \big( B_1(\lambda) + B_1(\lambda)^* \big) 
&\leq \big\| \til B_1(\lambda) \big\| (1+|\D|) R_\D(\lambda) .
\end{align*}
For $l=2$, we consider 
\[
\til B_2(\lambda) := (1+|\D|)^{-\frac12} \left( [\D,\rho]^* S \rho + \rho S [\D,\rho] \right) (1+|\D|)^{-\frac12} .
\]
Using the boundedness of $\left( [\D,\rho]^* S \rho + \rho S [\D,\rho] \right) (1+|\D|)^{-1}$ 
and \cref{lem:interpolation}, we see that $\til B_2(\lambda)$ is bounded. Since 
\[
B_2(\lambda) + B_2(\lambda)^* = (1+\lambda) R_\D(\lambda) (1+|\D|)^{\frac12} \til B_2(\lambda) (1+|\D|)^{\frac12} R_\D(\lambda) ,
\]
we obtain 
\begin{align*}
\pm \big( B_2(\lambda) + B_2(\lambda)^* \big) 
&\leq \big\| \til B_2(\lambda) \big\| (1+\lambda) (1+|\D|) R_\D(\lambda)^2 
\leq \big\| \til B_2(\lambda) \big\| (1+|\D|) R_\D(\lambda) ,
\end{align*}
where we have used that $\|(1+\lambda)R_\D(\lambda)\|\leq1$. 
For $l=3$, we consider 
\[
\til B_3(\lambda) := \D R_\D(\lambda)^{\frac12} (1+|\D|)^{-\frac12} \left( \rho S [\D,\rho] + [\D,\rho]^* S \rho \right) (1+|\D|)^{-\frac12} \D R_\D(\lambda)^{\frac12} .
\]
By a similar argument, along with the boundedness of $\D R_\D(\lambda)^{\frac12}$, we see that $\til B_3(\lambda)$ is bounded, and we obtain
\begin{align*}
\pm \big( B_3(\lambda) + B_3(\lambda)^* \big) 
&\leq \big\| \til B_3(\lambda) \big\| (1+|\D|) R_\D(\lambda) .
\end{align*}
Summarising, we have shown for $l=1,2,3$ that 
\begin{align*}
\pm \big( B_l(\lambda) + B_l^*(\lambda) \big) \leq C_l (1+|\D|) R_\D(\lambda) ,
\end{align*}
with the constants $C_l := \|\til B_l(\lambda)\|$. 
Inserting these inequalities 
into the definition of $I(\psi)$, we obtain
\begin{align*}
\pm I(\psi) \geq - \sum_{l=1}^3 C_l \Big\la \psi \Bigmvert \int_0^\infty \int_0^\infty (\lambda\mu)^{-\frac12} \Big( & (1+\mu) R_S(\mu) (1+|\D|) R_\D(\lambda) R_S(\mu) \\
&+ S R_S(\mu) (1+|\D|) R_\D(\lambda) S R_S(\mu) \Big) \psi d\lambda d\mu \Big\ra .
\end{align*}
By \cref{lem:integral_formula}, the integral over $\lambda$ converges strongly, and for the strong limit we have the norm bound 
\[
\left\| \int_0^\infty \lambda^{-\frac12} (1+|\D|) R_\D(\lambda) d\lambda \right\| \leq \big\| (1+|\D|) (1+\D^2)^{-\frac12} \big\| \leq 2 .
\]
Hence we obtain:
\begin{align*}
\pm I(\psi) &\geq - \sum_{l=1}^3 2 C_l \Big\la \psi \Bigmvert \int_0^\infty \mu^{-\frac12} \Big( (1+\mu) R_S(\mu)^2 + S^2 R_S(\mu)^2 \Big) \psi d\mu \Big\ra \\
&= - \sum_{l=1}^3 2 C_l \Big\la \psi \Bigmvert \int_0^\infty \mu^{-\frac12} R_S(\mu) \psi d\mu \Big\ra = - \sum_{l=1}^3 2 \pi C_l \big\la \psi \bigmvert (1+S^2)^{-\frac12} \psi \big\ra \\
&\geq - \sum_{l=1}^3 2 \pi C_l \la \psi | \psi \ra . 
\end{align*}
Thus we have proven the first statement with $C := \sum_{l=1}^3 2 \pi C_l$. 
The second statement follows immediately by observing that the operators $\til B_l(\lambda)$ (and hence the constants $C_l$) are linear in $S$. 
\end{proof}

\begin{lem}
\label{lem:integral_M}
The operator 
\[
M := \sum_{m=1}^4 \int_0^\infty \!\!\!\! \int_0^\infty \lambda^{-\frac12} \mu^{-\frac12} M_m(\lambda,\mu)^* \rho^2 M_m(\lambda,\mu) d\mu .
\]
is well-defined and bounded, and $\|M\| \leq 4\pi^2$. 
\end{lem}
\begin{proof}
For each $m=1,\ldots,4$, we have 
\[
\big\| M_m(\lambda,\mu) \big\| \leq \frac{1}{\sqrt{1+\lambda}\sqrt{1+\mu}} .
\]
Since $\|\rho\|\leq1$, we see that $\big\| M_m(\lambda,\mu)^* \rho^2 M_m(\lambda,\mu) \big\|$ is bounded by $(1+\lambda)^{-1}(1+\mu)^{-1}$, and therefore 
\begin{align*}
\|M\| 
&\leq \sum_{m=1}^4 \int_0^\infty \!\!\!\! \int_0^\infty \lambda^{-\frac12} \mu^{-\frac12} \big\| M_m(\lambda,\mu)^* \rho^2 M_m(\lambda,\mu) \big\| d\mu \\ 
&\leq 4 \int_0^\infty \!\!\!\! \int_0^\infty \lambda^{-\frac12} \mu^{-\frac12} (1+\lambda)^{-1}(1+\mu)^{-1} d\mu 
= 4\pi^2 .
\qedhere
\end{align*}
\end{proof}

\begin{prop}
\label{prop:local_positivity}
Let $0<\kappa<2$, and let $\chi=\bar\chi\in C_c^\infty(M)$ be such that $\rho|_{\supp\chi}=1$. 
Suppose there exists a constant $c>0$ such that for all $\psi\in\Dom\D$ we have 
\begin{align}
\label{eq:unbdd_pos_local}
\la\D \rho\psi|S\rho\psi\ra + \la S\rho\psi|\D \rho\psi\ra \geq - c \la \rho\psi|\rho\psi\ra .
\end{align}
Then there exists an $\alpha>0$ such that the operator $\chi [ F_\D , F_{\alpha S} ]_\pm \chi + \kappa \chi^2$ is positive modulo compact operators:
\[
\chi [F_\D , F_{\alpha S}]_\pm \chi \gtrsim - \kappa \chi^2 .
\]
\end{prop}
\begin{proof}
Let $\psi\in L^2(M,\bE)$. Using that $\rho^2\chi = \chi$, we can insert the result of \cref{lem:positivity_local_expr} into \cref{eq:comm_integral} to obtain
\begin{align*}
\big\la \psi & \bigmvert \chi [F_\D , F_S]_\pm \chi \psi \big\ra \\
&= \frac{1}{\pi^2} \int_0^\infty \int_0^\infty \lambda^{-\frac12} \mu^{-\frac12} \bigg( 2 \Re \big\la \chi\psi \bigmvert K(\lambda,\mu) \chi\psi \big\ra \\
&\quad- 2 \sum_{l=1}^3 \Re \big\la B_l(\lambda) \sqrt{1+\mu} R_S(\mu) \chi\psi \bigmvert \sqrt{1+\mu} R_S(\mu) \chi\psi \big\ra \\
&\quad- 2 \sum_{l=1}^3 \Re \big\la B_l(\lambda) S R_S(\mu) \chi\psi \bigmvert S R_S(\mu) \chi\psi \big\ra 
+ \sum_{m=1}^4 Q\big(\rho M_m(\lambda,\mu)\chi\psi\big) \bigg) d\lambda d\mu \\
&= \frac{1}{\pi^2} 2 \Re \big\la \chi\psi \bigmvert K \chi\psi \big\ra 
- \frac{1}{\pi^2} I(\chi\psi) 
+ \frac{1}{\pi^2} \sum_{m=1}^4 \int_0^\infty \int_0^\infty \lambda^{-\frac12} \mu^{-\frac12} Q\big(\rho M_m(\lambda,\mu)\chi\psi\big) d\lambda d\mu , 
\end{align*}
where we have inserted the definitions of $K$ and $I$. 
From \cref{lem:lower_bound_integral}, we have the lower bound 
\begin{align*}
- \frac{1}{\pi^2} I(\chi\psi) &\geq -\frac{C}{\pi^2} \big\la \chi \psi \bigmvert \chi \psi \big\ra .
\end{align*}
By \cref{eq:unbdd_pos_local}, there exists a constant $c>0$ such that 
\[
Q\big(\rho M_m(\lambda,\mu)\chi \psi\big) \geq - c \big\la \rho M_m(\lambda,\mu)\chi \psi \bigmvert \rho M_m(\lambda,\mu)\chi \psi \big\ra . 
\]
Using \cref{lem:integral_M} we obtain that 
\[
\frac{1}{\pi^2} \sum_{m=1}^4 \int_0^\infty \!\!\!\! \int_0^\infty \lambda^{-\frac12} \mu^{-\frac12} Q\big(fM_m(\lambda,\mu)\chi \psi\big) d\lambda d\mu 
\geq - \frac{1}{\pi^2} c \|M\| \big\la \chi \psi \bigmvert \chi \psi \big\ra \\
\geq - 4 c \big\la \chi \psi \bigmvert \chi \psi \big\ra .
\]
Hence we have shown that 
\[
\big\la \psi \bigmvert \chi [F_\D , F_S]_\pm \chi \psi \big\ra 
\geq \frac{1}{\pi^2} 2 \Re \big\la \chi \psi \bigmvert K \chi \psi \big\ra - \left( \frac{C}{\pi^2}  + 4 c \right) \big\la \chi \psi \bigmvert \chi \psi \big\ra .
\]
Since this holds for any $\psi$, we have the operator inequality
\[
\chi [F_\D , F_S]_\pm \chi \geq \frac{1}{\pi^2} \chi (K+K^*) \chi - \left( \frac{C}{\pi^2}  + 4 c \right) \chi^2 .
\]
Since $K$ is compact by \cref{lem:integral_K}, we have therefore shown that 
\[
\chi [F_\D , F_S]_\pm \chi \gtrsim - \left( \frac{C}{\pi^2}  + 4 c \right) \chi^2 .
\]
Finally, if we replace $S$ by $\alpha S$ for some $\alpha>0$, then $c$ and $C$ are replaced by $\alpha c$ and $\alpha C$. 
Thus, by choosing $\alpha$ small enough, we can ensure that $\alpha ( C\pi^{-2}  + 4 c ) < \kappa < 2$. 
\end{proof}

\section{The internal Kasparov product}
\label{sec:Kasp_prod}

In this section we will show that we can construct the Kasparov product of a vertical and a horizontal operator on a submersion. 
The proof is obtained by checking the connection and positivity conditions in the following well-known theorem by Connes and Skandalis. We cite below a slightly more general version of their theorem, as described in the comments following \cite[Definition 18.4.1]{Blackadar98}. 
For convenience, let us first introduce some notation. Given a Hilbert $B$-module $E_1$ and a Hilbert $C$-module $E_2$ with a $*$-homomorphism $B\to\End_C(E_2)$, we consider the internal tensor product $E := E_1\hot_BE_2$. For any $\psi\in E_1$, we define the operator $T_\psi \colon E_2 \to E$ as $T_\psi \eta = \psi\hot\eta$ for any $\eta\in E_2$. The operator $T_\psi$ is adjointable, and its adjoint $T_\psi^*\colon E\to E_2$ is given by $T_\psi^* (\xi\otimes\eta) = \la\psi|\xi\ra\cdot\eta$. 
Furthermore, we also introduce the operator $\til T_\psi$ on the Hilbert $C$-module $E\oplus E_2$ given by 
\[
\til T_\psi := \mattwo{0}{T_\psi}{T_\psi^*}{0} . 
\]

\begin{thm}[{\cite[Theorem A.3]{CS84} \& \cite[Theorem 18.4.3]{Blackadar98}}]
\label{thm:Connes-Skandalis}
Consider $C^*$-algebras $A,B,C$, where $A$ is separable and $B,C$ are $\sigma$-unital. 
Let $(A,{}_{\phi_1}{E_1}_B,F_1)$ and $(B,{}_{\phi_2}{E_2}_C,F_2)$ be Kasparov modules, and consider the Hilbert $C$-module $E := E_1\hot_B E_2$ and the $*$-homomorphism $\phi := \phi_1\hot1 \colon A \to \End_C(E)$. 
Suppose that $(A,{}_\phi E_C,F)$ is a Kasparov module such that the following two conditions hold:
\begin{description}
\item[Connection condition:] for any $\psi\in E_1$, 
the graded commutator $[F\oplus F_2,\til T_\psi]_\pm$ is compact on $E\oplus E_2$; 
\item[Positivity condition:] there exists a $0\leq\kappa<2$ such that for all $a\in A$ we have that $\phi(a) [F_1\hot1,F]_\pm \phi(a^*) + \kappa\phi(aa^*)$ is positive modulo compacts on $E$. 
\end{description}
Then $(A,{}_\phi E_C,F)$ represents the Kasparov product of $(A,{}_{\phi_1}{E_1}_B,F_1)$ and $(B,{}_{\phi_2}{E_2}_C,F_2)$:
\[
[F] = [F_1] \hot_B [F_2] \in \KK(A,C) .
\]
Moreover, an operator $F$ with the above properties always exists and is unique up to operator homotopy. 
\end{thm}

Now let $\pi\colon M\to B$ be a submersion of Riemannian manifolds, and let $\D_V$ on $\bE_V\to M$ be as in \cref{ass:vert_op}. We assume that the bundle $\bE_V = \bE_V^+\oplus\bE_V^-$ is $\Z_2$-graded, and that $\D_V$ is an odd operator. 
Consider furthermore an odd first-order symmetric \emph{elliptic} differential operator $\D_B$ on a $\Z_2$-graded hermitian vector bundle $\bE_B$ over the base manifold $B$, which yields a $\K$-homology class $[\D_B] \in \KK(C_0(B),\C)$. 

We consider the `horizontal' bundle $\bE_H := \pi^* \bE_B$ on $M$, and consider the tensor product 
\[
\Gamma_0(B,E_\bullet) \hot_{C_0(B)} L^2(B,\bE_B) \simeq L^2(M,\bE_V\hot\bE_H) .
\]
The operator $\D_V$ gives a symmetric operator $\D_V\hot1$ on $\Gamma_0(B,E_\bullet) \hot_{C_0(B)} L^2(B,\bE_B)$. 
Consider a given hermitian connection $\nabla$ on the $C_c^\infty(B)$-module $\Gamma_c^\infty(M,\bE_V)$, i.e.\ a map $\nabla\colon \Gamma_c^\infty(M,\bE_V) \to \Gamma_c^\infty(M,\bE_V) \otimes_{C_c^\infty(B)} \Gamma_c^\infty(B,T^*B)$ 
satisfying the Leibniz rule $\nabla(\psi \pi^*(f)) = \nabla(\psi) f + \psi\otimes df$ for $\psi \in \Gamma_c^\infty(M,\bE_V)$ and $f\in C_c^\infty(B)$, 
and the hermitian property $(\nabla_X\psi_1|\psi_2) + (\psi_1|\nabla_X\psi_2) = X (\psi_1|\psi_2)$ for all $\psi_1,\psi_2 \in \Gamma_c^\infty(M,\bE_V)$ and $X\in\Gamma^\infty(TB)$. 
We then define the operator $1\hot_\nabla\D_B$ on the initial domain $\Gamma_c^\infty(M,\bE_V) \hot_{C_c^\infty(B)} \Gamma_c^\infty(B,\bE_B)$ in $\Gamma_0(B,E_\bullet) \hot_{C_0(B)} L^2(B,\bE_B)$ by 
\begin{align}
\label{eq:lift_D_B}
(1\hot_\nabla\D_B)(\psi\hot\eta) := \psi\hot\D_B\eta + \sum_{k=1}^{\dim B} \nabla_{e_k} \psi \hot \sigma_B(e_k) \eta ,
\end{align}
where $\{e_k\}$ is an orthonormal frame of $TB$ and $\sigma_B$ is the principal symbol of $\D_B$. Since the connection $\nabla$ is hermitian, the operator $1\hot_\nabla\D_B$ is again symmetric. 
We consider the \emph{tensor sum} 
\[
\D := \D_V\hot1 + 1\hot_{\nabla}\D_B ,
\]
which we view as a symmetric first-order differential operator on $L^2(M,\bE_V\hot\bE_H)$ with initial domain $\Gamma_c^\infty(M,\bE_V\hot\bE_H)$. 

\begin{lem}
The tensor sum $\D$ is elliptic. 
\end{lem}
\begin{proof}
Since $\D_V$ and $\D_B$ are odd operators, the principal symbols $\sigma_V$ and $\sigma_B$ of $\D_V\hot1$ and $1\hot_\nabla\D_B$ (respectively) anti-commute. Hence the square of the principal symbol $\sigma_\D$ of $\D$ is given by $\sigma_\D(x,\xi)^2 = \sigma_V(x,\xi_V)^2 + \sigma_B(x,\xi_H)^2$ for any $x\in M$ and $\xi = \xi_V\oplus\xi_H \in (TM)_x = (T_VM)_x \oplus (T_HM)_x$. 
Since $\D_V\hot1$ and $1\hot_\nabla\D_B$ are symmetric, we know that both $\sigma_V(x,\xi_V)^2$ and $\sigma_B(x,\xi_H)^2$ are positive. 
Since $\D_V$ is vertically elliptic and $\D_B$ is elliptic on the base, we then see that $\sigma_\D(x,\xi)^2$ is invertible for any $0\neq\xi\in(TM)_x$. 
Thus $\D$ is elliptic. 
\end{proof}

It follows that $(C_c^\infty(M) , L^2(M,\bE_V\hot\bE_H) , \D)$ is a half-closed module, and we obtain a class $[\D] \in \KK(C_0(M),\C)$. 
Our aim is to prove that the operator $\D$ represents the Kasparov product of $\D_V$ and $\D_B$:
\[
[\D] = [\D_V] \hot_{C_0(B)} [\D_B] \in \KK(C_0(M),\C) .
\]
We consider a locally finite cover $\{U_j\}$ of $M$ with a corresponding partition of unity $\{\chi_j^2\}$ and functions $\{\phi_j\}$, satisfying the same conditions as in \cref{ass:Higson}. We represent the $\KK$-class of $\D_V$ by the localised representative as constructed in \cref{defn:Higson}: 
\[
\til F_V(\alpha) := \sum_j \chi_j F_{\alpha_j\D_{V,j}} \chi_j ,
\]
for some sequence $\{\alpha_j\}$ of strictly positive numbers, where $\D_{V,j} := \phi_j\D_V\phi_j$. 
The classes of $\D$ and $\D_B$ are represented simply by their bounded transforms $F_\D$ and $F_B := F_{\D_B}$. 
Thus we aim to prove that 
\[
[F_\D] = [\til F_V(\alpha)] \hot_{C_0(B)} [F_B] \in \KK(C_0(M),\C) .
\]
The remainder of this section is devoted to the proof of this equality in the following setting:
\begin{assumption}
\label{ass:Kasp_prod}
Let $\pi\colon M\to B$ be a submersion of Riemannian manifolds. 
Let $\D_V$ be an odd vertically elliptic symmetric first-order differential operator on a $\Z_2$-graded hermitian vector bundle $\bE_V\to M$, such that for each $b\in B$, the subspace $\ev_b(\Dom\D_V^*) \subset \Dom(\D_V)_b^*$ is a \emph{core} for $(\D_V)_b^*$. 
Let $\D_B$ be an odd symmetric elliptic first-order differential operator on a $\Z_2$-graded hermitian vector bundle $\bE_B\to B$. 
Let $\nabla$ be a hermitian connection on the $C_c^\infty(B)$-module $\Gamma_c^\infty(M,\bE_V)$. 
\end{assumption}

\begin{prop}[Connection condition]
\label{prop:connection}
For $\psi\in\Gamma_c^\infty(M,\bE_V)$, the operator 
$[ F_{\D\oplus\D_B} , \til T_\psi ]_\pm$ 
is compact on the Hilbert space $L^2(M,\bE_V\hot\bE_H) \oplus L^2(B,\bE_B)$. 
\end{prop}
\begin{proof}
Using \cref{lem:integral_formula}, we have a strongly convergent integral 
\begin{align*}
\big[ F_{\D\oplus\D_B} , \til T_\psi \big]_\pm 
&= - \frac1\pi \int_0^\infty \lambda^{-1/2} \big[ (\D\oplus\D_B) R_{\D\oplus\D_B}(\lambda) , \til T_\psi \big]_\pm d\lambda .
\end{align*}
It is a standard computation to check that $\til T_\psi\in\Lip(\D\oplus\D_B)$ (in our context, see e.g.\ the proof of \cite[Theorem 22]{KS17bpre}). 
By \cref{lem:R_S_cpt_order}, we then know that the operator $\big[ (\D\oplus\D_B) R_{\D\oplus\D_B}(\lambda) , \til T_\psi \big]_\pm$ is compact and of order $\mO(\lambda^{-1})$. Hence the above integral is in fact norm-convergent, and therefore $[ F_{\D\oplus\D_B} , \til T_\psi ]_\pm$ is compact. 
\end{proof}

\begin{lem}
\label{lem:unbdd_positivity}
For each compact subset $K\subset M$, there exists a constant $c_K>0$ such that for all $\psi\in \Gamma_c^\infty(M,\bE_V\hot\bE_H)$ with $\supp(\psi)\subset K$ we have 
\[
\la\D\psi|(\D_V\hot1)\psi\ra + \la(\D_V\hot1)\psi|\D\psi\ra \geq - c_K \la\psi|\psi\ra .
\]
\end{lem}
\begin{proof}
The proof follows the results of \cite{KS17bpre}; here we only give a brief sketch. First, as in \cite[Lemma 16]{KS17bpre}, one shows that the anti-commutator $[\D_V\hot1,1\hot_{\nabla}\D_B]_\pm$ is a vertical first-order differential operator. 
Since $\D_V$ is vertically elliptic, an application of G\aa rding's inequality then shows that, given the compact subset $K\subset M$, there exists a constant $c_K'>0$ such that for every $\psi\in\Gamma_c^\infty(M,\bE_V\hot\bE_H)$ with $\supp(\psi)\subset K$ we have $\|[\D_V\hot1,1\hot_{\nabla}\D_B]_\pm\psi\| \leq c_K' \|\psi\|_{\D_V\otimes1}$ (see \cite[Lemma 17]{KS17bpre}). 
Following the argument given in the proof of \cite[Theorem 22]{KS17bpre}, we find that the statement holds with the constant $c_K = \frac12(1+c_K')$. 
\end{proof}

\begin{prop}[Positivity condition]
\label{prop:positivity}
Let $0<\kappa<2$. 
There exists a sequence $\{\alpha_j\}_{j\in\N}$ of positive real numbers such that for any $f \in C_0(M)$, 
we have that $f [F_\D , \til F_V(\alpha)\hot1]_\pm \bar f + \kappa f\bar f$ is positive modulo compact operators. 
\end{prop}
\begin{proof}
Since $C_c^\infty(M)$ is dense in $C_0(M)$, it suffices to prove the statement for $f\in C_c^\infty(M)$. 
In this case, we know that 
$\sum_j f \chi_j$ 
is a finite sum. 
Using that $[F_\D,\chi_j] \sim 0$ by \cref{thm:Hilsum} and $\chi_j (F_\D - F_{\D_j}) \sim 0$ by \cref{lem:local_comparison} (where $\D_j := \phi_j \D \phi_j$), we have 
\begin{align*}
f [ F_\D , \til F_V(\alpha)\hot1 ]_\pm \bar f 
&= \sum_j f [ F_\D , \chi_j (F_{\alpha_j\D_{V,j}}\hot1) \chi_j ]_\pm \bar f 
\sim \sum_j f \chi_j [ F_\D , F_{\alpha_j\D_{V,j}}\hot1 ]_\pm \chi_j \bar f \\
&\sim \sum_j f \chi_j [ F_{\D_j} , F_{\alpha_j\D_{V,j}}\hot1 ]_\pm \chi_j \bar f .
\end{align*}

We would like to apply the results from \cref{sec:local_positivity} to the operators $\chi_j [ F_{\D_j} , F_{\alpha_j\D_{V,j}}\hot1 ]_\pm \chi_j$. We pick a function $\rho_k\in C_c^\infty(M,[0,1])$ such that $\phi_k|_{\supp\rho_k} = 1$ and $\rho_k|_{\supp\chi_k}=1$. 
We need to check that \cref{ass:local_positivity} is satisfied by $\D_j$, $S_j = \D_{V,j}\hot1$, and $\rho_j$. 
Consider the domain $\Dom\D_j = \Dom(\D\phi_j^2) = \{ \psi\in L^2(M,\bE_V\hot\bE_H) : \phi_j^2\psi\in\Dom\D \}$. Since $\phi_j^2\psi$ is compactly supported, an application of G\aa rding's inequality shows that $\phi_j^2\psi \subset \Dom\D_V$ (see \cite[Lemma 21]{KS17bpre}). Therefore we have the domain inclusion $\Dom\D_j = \Dom(\D\phi_j^2) \subset \Dom(\D_V\phi_j^2) = \Dom\D_{V,j}$. 
Moreover, by construction, $\D_j$ is elliptic on (a neighbourhood of) the support of $\rho_j$. 

From \cref{lem:unbdd_positivity} we know that for each $j$ there exists a constant $c_j>0$ such that for all $\psi\in \Gamma_c^\infty(M,\bE_V\hot\bE_H)$ we have 
\[
\la\D_j\chi_j\psi|(\D_{V,j}\hot1)\chi_j\psi\ra + \la(\D_{V,j}\hot1)\chi_j\psi|\D_j\chi_j\psi\ra \geq - c_j \la\chi_j\psi|\chi_j\psi\ra ,
\]
where we used that $\D\chi_j = \D_j\chi_j$ and $\D_V\chi_j = \D_{V,j}\chi_j$. 
So by \cref{prop:local_positivity} there exists a sequence $\{\alpha_j\}_{j\in\N} \subset (0,\infty)$ such that $\chi_j [ F_{\D_j} , F_{\alpha_j\D_{V,j}}\hot1 ]_\pm \chi_j + \kappa \chi_j^2$ is positive modulo compact operators. Hence we have
\[
f [ F_\D , \til F_V(\alpha)\hot1 ]_\pm \bar f 
\sim \sum_j f \chi_j [ F_{\D_j} , F_{\alpha_j\D_{V,j}}\hot1 ]_\pm \chi_j \bar f 
\gtrsim - \sum_j \kappa f \chi_j^2 \bar f 
= - \kappa f\bar f .
\]
Thus we obtain that $f [ F_\D , \til F_{\D_1}(\alpha)\hot1 ]_\pm \bar f  + \kappa f\bar f$ is positive modulo compact operators. 
\end{proof}

\begin{thm}
\label{thm:submersion_Kasp_prod}
Consider the setting of \cref{ass:Kasp_prod}. 
Then the tensor sum $\D$ represents the Kasparov product of $\D_V$ and $\D_B$:
\[
[\D] = [\D_V] \hot_{C_0(B)} [\D_B] \in \KK(C_0(M),\C) .
\]
\end{thm}
\begin{proof}
We have $[\D] = [F_\D]$, $[\D_B] = [F_B]$, and $[\D_V] = [\til F_V(\alpha)]$ (by \cref{thm:local_rep_KK}). 
Using \cref{prop:connection,prop:positivity}, the statement then follows from \cref{thm:Connes-Skandalis}. 
\end{proof}

\section{Factorisation of the Dirac operator on a submersion}
\label{sec:factorisation}

In this section we consider a submersion $\pi\colon M\to B$ of smooth even-dimensional Riemannian spin$^c$ manifolds $M$ and $B$. 
We recall that $T_VM = \Ker d\pi$ denotes the vertical tangent bundle of $M$, and the horizontal tangent bundle is then given by the orthogonal complement $T_HM := (T_VM)^\perp$. 
We will assume furthermore that the submersion is \emph{Riemannian}, which means that $d\pi(x) \colon (T_HM)_x \to (TB)_{\pi(x)}$ is an isometry for all $x\in M$. 
We aim to prove that the Dirac operator $\D_M$ on the total manifold $M$ can be factorised in unbounded $\KK$-theory in terms of a vertical operator $\D_V$ and the Dirac operator $\D_B$ on the base manifold $B$, up to an explicit curvature term. 
We closely follow the work of Kaad and Van Suijlekom \cite{KS18,KS17bpre}, who proved this factorisation result for a \emph{proper} submersion (i.e.\ when each fibre $M_b = \pi^{-1}(b)$ is compact). 

Let $\bS_M\to M$ be the smooth $\Z_2$-graded spinor bundle over $M$. Since $M$ is spin$^c$, the Clifford multiplication yields an even isomorphism 
\[
c_M \colon \Cliff(M) \to \End_{C^\infty(M)}\big(\Gamma^\infty(M,\bS_M)\big) ,
\]
where $\Cliff(M) = \Gamma^\infty(M,\Cliff(TM))$ denotes the Clifford algebra over $M$. The Levi-Civita connection can be lifted to an even hermitian Clifford connection $\nabla^{\bS_M}$ on $\Gamma^\infty(M,\bS_M)$. The Dirac operator is then defined by
\begin{align}
\label{eq:Dirac}
\D_M := c_M \circ \nabla^{\bS_M} = \sum_{j=1}^{\dim M} c_M(e_j) \nabla^{\bS_M}_{e_j} ,
\end{align}
where $\{e_j\}$ is a local orthonormal frame for $TM$. Similarly, we also have a spinor bundle $\bS_B\to B$ and a Dirac operator $\D_B := c_B \circ \nabla^{\bS_B}$. Since both Dirac operators are elliptic and symmetric, we obtain a half-closed $C_0(M)$-$\C$-module $(C_c^\infty(M) , L^2(M,\bS_M), \D_M )$ and a half-closed $C_0(B)$-$\C$-module $(C_c^\infty(B) , L^2(B,\bS_B) ,\D_B )$. 

Let $\Cliff_V(M) := \Gamma^\infty(M,\Cliff(T_VM))$ and $\Cliff_H(M) := \Gamma^\infty(M,\Cliff(T_HM))$ denote the Clifford algebras of vertical and horizontal vector fields, respectively. 
We pull back the spinor bundle over $B$ to a \emph{horizontal spinor bundle} $\bS_H := \pi^*\bS_B$ over $M$, which is equipped with the Clifford multiplication $c_H$ by the horizontal vector fields and with a hermitian Clifford connection $\nabla^{\bS_H}$. 
We then define the \emph{vertical spinor bundle} 
\[
\bS_V := \bS_H^* \otimes_{\Cliff(T_HM)} \bS_M 
\]
which is equipped with the Clifford multiplication by vertical vector fields $c_V\colon \Cliff_V(M) \to \Gamma^\infty(M,\End\bS_V)$ and with a hermitian Clifford connection $\nabla^{\bS_V}$. 
We note that we have a natural isomorphism $\bS_M \simeq \bS_V \hot \bS_H$. 
For more details on these constructions and explicit formulae, we refer to \cite[\S3]{KS18}. 
The \emph{vertical Dirac operator} $\D_V$ is then defined as
\begin{align}
\label{eq:vertical_Dirac}
\D_V := c_V \circ \nabla^{\bS_V} = \sum_{j=1}^{\dim M - \dim B} c_V(e_j) \nabla^{\bS_V}_{e_j} ,
\end{align}
where $\{e_j\}$ is a local orthonormal frame of $T_VM$. As in \cite[Lemma 12 \& Proposition 13]{KS17bpre}, we see that $\D_V$ is an odd vertically elliptic symmetric first-order differential operator. 
We will view $\D_V$ as an odd symmetric operator on the $\Z_2$-graded Hilbert $C_0(B)$-module $\Gamma_0(B,E_\bullet)$, where $E_\bullet$ denotes the bundle of Hilbert spaces $E_b := L^2(M_b,\bS_V|_{M_b})$ for $b\in B$. 
If (the closure of) $\D_V$ is regular, we obtain a half-closed $C_0(M)$-$C_0(B)$-module $(C_c^\infty(M) , \Gamma_0(B,E_\bullet) , \D_V)$. 
We then also obtain an odd regular symmetric operator $\D_V\hot 1$ on the interior tensor product $\Gamma_0(B,E_\bullet) \hot_{C_0(B)} L^2(B,\bS_B)$, which is isomorphic to the Hilbert space $L^2(M,\bS_M)$. 

We also want to lift the operator $\D_B$ to an operator on $\Gamma_0(B,E_\bullet) \hot_{C_0(B)} L^2(M,\bS_H)$, and for this purpose we need to choose a connection on $\Gamma_0(B,E_\bullet)$. 
As in \cite[Definition 17]{KS18}, we use the \emph{mean curvature} $k \in \Hom_{C^\infty(M)}\big(\Gamma^\infty(T_HM),C^\infty(M)\big)$ to define a hermitian connection $\nabla$ on $\Gamma_c^\infty(M,\bS_V)$ by
\begin{align}
\label{eq:vert_connection}
\nabla_X \psi &:= \nabla^{\bS_V}_{X_H} \psi + \frac12 k(X_H) \cdot \psi ,
\end{align}
for $\psi \in \Gamma_c^\infty(M,\bS_V)$ and for a vector field $X$ on $B$ with horizontal lift $X_H \in \Gamma^\infty(M,\pi^*TB)$. We then consider the operator $1\hot_\nabla\D_B$ on the initial domain $\Gamma_c^\infty(M,\bS_V) \hot_{C_c^\infty(B)} \Gamma_c^\infty(B,\bS_B)$ in $\Gamma_0(B,E_\bullet) \hot_{C_0(B)} L^2(B,\bS_B)$ defined by \cref{eq:lift_D_B}. 
As in \cite[Lemma 20]{KS18}, the operator $1\hot_\nabla\D_B$ is odd and symmetric. 

Consider the curvature form $\Omega \in \Gamma^\infty(M,T_H^*M\wedge T_H^*M\otimes T_V^*M)$ of the Riemannian submersion $\pi\colon M\to B$ given by 
\[
\Omega(X,Y,Z) := g_M([X,Y],Z) ,
\]
where $g_M$ is the Riemannian metric on $M$, $X,Y$ are horizontal vector fields, and $Z$ is a vertical vector field. This curvature form acts as an endomorphism on the bundle $\bS_M$ via the map $c\colon \Gamma^\infty(M,T_H^*M\wedge T_H^*M\otimes T_V^*M) \to \Gamma^\infty(M,\End\bS_M)$ given by
\[
c(\omega_1\wedge\omega_2\otimes\omega_3) := \big[ c_M(\omega_1^\sharp) , c_M(\omega_2^\sharp) \big] \cdot c_M(\omega_3^\sharp) ,
\]
where we have used the `musical' isomorphisms $\sharp \colon \Gamma^\infty(M,T_H^*M) \to \Gamma^\infty(M,T_HM)$ and $\sharp \colon \Gamma^\infty(M,T_V^*M) \to \Gamma^\infty(M,T_VM)$. We thus obtain a smooth endomorphism $c(\Omega)$. Since $M$ is non-compact, we note that $c(\Omega)$ is not necessarily globally bounded. 

Using \cref{thm:submersion_Kasp_prod}, we now obtain a generalisation of \cite[Theorem 22]{KS17bpre} to the case of a Riemannian submersion where the fibres are allowed to be non-compact and incomplete. 

\begin{thm}
\label{thm:Dirac_factorisation}
Let $\pi\colon M\to B$ be a Riemannian submersion of even-dimensional Riemannian spin$^c$ manifolds. 
Suppose that the vertical operator $\D_V$ defined in \cref{eq:vertical_Dirac} 
has the property that for each $b\in B$ the subspace $\ev_b(\Dom\D_V^*) \subset \Dom(\D_V)_b^*$ is a \emph{core} for $(\D_V)_b^*$. 
Then the half-closed module $( C_c^\infty(M) , L^2(M,\bS_M) , \D_M )$ is the unbounded Kasparov product of the half-closed module $( C_c^\infty(M) , \Gamma_0(B,E_\bullet) , \D_V )$ with the half-closed module $( C_c^\infty(B) , L^2(B,\bS_B) , \D_B )$ up to the curvature term $-\frac i8 c(\Omega)$. 
\end{thm}
\begin{proof}
The assumption on $\D_V$ ensures that we obtain from \cref{prop:vert_half-closed} a half-closed module $( C_c^\infty(M) , \Gamma_0(B,E_\bullet) , \D_V )$. 
By \cref{thm:submersion_Kasp_prod}, $\D$ represents the Kasparov product of $\D_V$ and $\D_B$. 
By \cite[Proposition 18]{KS17bpre}, $\D$ is unitarily equivalent to $\D_M - \frac i8 c(\Omega)$ under the unitary isomorphism $\Gamma_0(B,E_\bullet) \hot_{C_0(B)} L^2(B,\bS_B) \simeq L^2(M,\bS_M)$. 
Finally, we know from \cref{prop:loc_bdd} that $\D_M$ and $\D_M - \frac i8 c(\Omega)$ represent the same class in the $\K$-homology of $M$. 
\end{proof}

\begin{example}[Two-dimensional domains]
We return to the setting of \cref{sec:regular_examples}. 
Let $M$ be an open subset of $(-1,1)\times\R$. We have a natural map $\pi\colon M\to(-1,1)$ given by $\pi(x,y) := x$. We equip $(-1,1)$ with the standard flat metric and $M$ with a metric of the form
\[
g_M(x,y) = dx^2 + h(x,y) dy^2 ,
\]
where $h$ is a smooth, strictly positive function on $M$. We assume that for each $x\in\R$, $\pi^{-1}(x)$ is not empty, so that $\pi$ is a Riemannian submersion. 
The mean curvature $k\in C^\infty(M)$ is in this case explicitly given by $k = \frac12 h^{-1} \partial_x h$, and the curvature form $\Omega$ vanishes identically. 

We consider the vertical Dirac operator $\D_V$ on $C_c^\infty(M)$, the `horizontal' Dirac operator $\D_B$ on $C_c^\infty((-1,1))$, and the total Dirac operator $\D_M$ on $C_c^\infty(M,\C^2)$ given by 
\begin{align*}
\D_V &:= -i \sqrt{h}^{-1} \partial_y , & 
\D_B &:= -i\partial_x , &
\D_M &:= \mattwo{0}{\partial_x-i \sqrt{h}^{-1} \partial_y}{-\partial_x-i \sqrt{h}^{-1} \partial_y}{0} .
\end{align*}
We \emph{assume} that the vertical operator $\D_V$ has the property that for each $x\in(-1,1)$ the subspace $\ev_x(\Dom\D_V^*) \subset \Dom(\D_V)_x^*$ is a \emph{core} for $(\D_V)_x^*$. 
From \cref{prop:vert_half-closed} we obtain a half-closed $C_0(M)$-$C_0(-1,1)$-module $( C_c^\infty(M) , \Gamma_0((-1,1),E_\bullet) , \D_V )$, a half-closed $C_0(-1,1)$-$\C$-module $( C_0(-1,1) , L^2(-1,1) , \D_B )$, and a half-closed $C_0(M)$-$\C$-module $( C_c^\infty(M) , L^2(M,\C^2) , \D_M )$, representing classes $[\D_V] \in \KK_1(C_0(M),C_0(-1,1))$, $[\D_B] \in \KK_1(C_0(-1,1),\C)$, and $[\D_M] \in \KK(C_0(M),\C)$, respectively. 
By \cref{thm:Dirac_factorisation}, $\D_M$ represents the Kasparov product of $\D_V$ and $\D_B$: 
\[
[\D_M] = [\D_V] \otimes_{C_0(-1,1)} [\D_B] \in \KK(C_0(M),\C) .
\]
Note that, although we have considered the even-dimensional case throughout this article, here $\D_V$ and $\D_B$ are \emph{odd}-dimensional. In this case, the tensor sum describing the Kasparov product (cf.\ \cite{KL13}) is given by 
\[
\D := \mattwo{0}{\D_V\hot1 + i\hot_\nabla\D_B}{\D_V\hot1 - i\hot_\nabla\D_B}{0} ,
\]
where $\nabla$ is defined as in \cref{eq:vert_connection}. 
Since the curvature form $\Omega$ vanishes identically, the factorisation in unbounded $\KK$-theory is in this case \emph{exact}, in the sense that $\D_M$ is unitarily equivalent to the tensor sum $\D$. 
\end{example}

%\bibliographystyle{myamsalpha}
%\bibliography{short,bibliography}

\providecommand{\noopsort}[1]{}
\providecommand{\bysame}{\leavevmode\hbox to3em{\hrulefill}\thinspace}
\providecommand{\MR}{\relax\ifhmode\unskip\space\fi MR }
\providecommand{\MRhref}[2]{%
  \href{http://www.ams.org/mathscinet-getitem?mr=#1}{#2}
}
\providecommand{\href}[2]{#2}
\providecommand{\doilinktitle}[2]{#1}
\providecommand{\doilinkjournal}[2]{\href{https://doi.org/#2}{#1}}
\providecommand{\doilinkvynp}[2]{\href{https://doi.org/#2}{#1}}
\providecommand{\eprint}[2]{#1:\href{https://arxiv.org/abs/#2}{#2}}

\end{document}